\documentclass[12pt]{amsart}
\usepackage{amsmath}
\usepackage{amssymb}
\usepackage{latexsym}
\usepackage[all]{xy}
\usepackage{caption}
\usepackage{graphicx,psfrag}
\usepackage{multirow}
\usepackage{color}
\usepackage{epstopdf}

\newtheorem{tpcl}{Tpcl}[section]
\newtheorem{Lemma}[tpcl]{Lemma}
\newtheorem{Theorem}[tpcl]{Theorem}
\newtheorem{Proposition}[tpcl]{Proposition}

\newtheorem{Definition}[tpcl]{Definition}
\newtheorem{Example}[tpcl]{Example}
\newtheorem{Remark}[tpcl]{Remark}

\begin{document}

\title[Immersed surfaces of minimal genus]{Immersed surfaces of minimal genus representing mod 2 homology classes of non-hyperbolic geometric closed 3-manifolds}

\author{Xiaoming Du}
\address{School of Mathematics, South China University of Technology,
Guangzhou, 510640, P. R. China}
\email{scxmdu@scut.edu.cn}
\thanks{This authors are partially supported by NSFC (Grant No.~12671085 and No.~12331003).}

\author{Xuezhi Zhao}
\address{School of Mathematical Sciences, Capital Normal University,
Beijing, 100048, P. R. China}
\email{zhaoxve@mail.cnu.edu.cn}

\date{}

%  Math Subject Classifications
\renewcommand{\subjclassname}{\textup{2020\ }Mathematics Subject Classification}

\maketitle

\begin{abstract}
Let $M$ be a closed orientable non-hyperbolic geometric 3-manifold. Some elements of $H_2(M;\mathbb{Z}_2)$ (here $\mathbb{Z}_2$ denotes $\mathbb{Z}/2\mathbb{Z}$) can only be represented by non-orientable closed surfaces, and an immersed representing surface can be much simpler than an embedded one. In this paper, for every element of $H_2(M;\mathbb{Z}_2)$, we determine the minimal complexity of an immersed representing surface and construct such a surface. As a by-product, we obtain an example in which this simplest representing surface is connected but not $\pi_1$-injective.
\end{abstract}

\section{Introduction}\label{sec-intro}

The realization of homology classes by surfaces is a classical problem in low-dimensional topology. Let $M$ be a connected, closed, orientable, irreducible 3-manifold and let $\alpha$ be a non-trivial element of the second homology group. We want to find the simplest surface representing $\alpha$. Apart from a few low genus cases, the complexity of a surface is measured by its negative Euler characteristic.

When the coefficient ring is $\mathbb{Z}_2$ (throughout this paper, $\mathbb{Z}_2$ means $\mathbb{Z}/2\mathbb{Z}$), the elements of $H_2(M;\mathbb{Z}_2)$ outside the image of $H_2(M;\mathbb{Z})\otimes \mathbb{Z}_2$ can only be represented by non-orientable surfaces. A representing surface may be disconnected: if $f_1:F_1\to M$ and $f_2:F_2\to M$ represent $\alpha_1$ and $\alpha_2$ respectively, then the disjoint union $f:F_1 \sqcup F_2 \to M$ represents $\alpha_1+\alpha_2$. For a surface with several components we define its complexity to be the sum of the complexities of its components. The connected sum $F_1 \# F_2$ is not the simplest representing surface, since $-\chi(F_1\#F_2) = -\chi(F_1) - \chi(F_2)+2> -\chi(F_1) - \chi(F_2)$.

For the simplest representing surface, we can consider both embedded ones and immersed ones.

\begin{Definition}\label{def-norms}
We define the embedded norm and the immersed norm for classes in the second homology group as follows:
$$
\begin{array}{l}
\displaystyle \|\alpha\|^{\mathrm{em}} = \min_{F \hookrightarrow M,\ [F] = \alpha}\left\{\sum_{F_i \subset F} \max \{0,-\chi(F_i)\}\right\} \qquad (\alpha \in H_2(M;\mathbb{Z})), \\
\displaystyle \|\alpha\|^{\mathrm{em}}_{\mathbb{Z}_2} = \min_{F \hookrightarrow M,\ [F] = \alpha}\left\{\sum_{F_i \subset F} \max \{0,-\chi(F_i)\}\right\} \qquad (\alpha \in H_2(M;\mathbb{Z}_2)), \\
\displaystyle \|\alpha\|^{\mathrm{im}} = \min_{F \;\looparrowright \, M,\ [F] = \alpha}\left\{\sum_{F_i \subset F} \max \{0,-\chi(F_i)\}\right\} \qquad (\alpha \in H_2(M;\mathbb{Z})),\\
\displaystyle \|\alpha\|^{\mathrm{im}}_{\mathbb{Z}_2} = \min_{F \;\looparrowright \, M,\ [F] = \alpha}\left\{\sum_{F_i \subset F} \max \{0,-\chi(F_i)\}\right\} \qquad (\alpha \in H_2(M;\mathbb{Z}_2)),
\end{array}
$$
where the minimum is taken over all closed surfaces $F$ representing $\alpha$, embedded in the first two cases and immersed in the last two, and $F_i$ runs over the connected components of $F$.
\end{Definition}

By definition, components of a surface that are spheres $S^2$, projective planes $\mathbb{R}P^2$, tori $T^2$, or Klein bottles $K$ do not contribute to the above norms. For $\alpha \in H_2(M;\mathbb{Z})$, the norm $\|\alpha\|^{\mathrm{em}}$ is the Thurston norm, and by a result of Gabai \cite{Gabai} we have $\|\alpha\|^{\mathrm{em}} = \|\alpha\|^{\mathrm{im}}$. For $\alpha \in H_2(M;\mathbb{Z}_2)$, the norm $\|\alpha\|^{\mathrm{em}}_{\mathbb{Z}_2}$ is the $\mathbb{Z}_2$-Thurston norm introduced in \cite{JRT}, and $\|\alpha\|^{\mathrm{im}}_{\mathbb{Z}_2}$ is its immersed counterpart. Here the two norms may differ: for $M=L(2k,1)$ with $k \geqslant 3$ we have $H_2(M;\mathbb{Z}_2)=\mathbb{Z}_2$, and for $\alpha \neq 0 \in H_2(M;\mathbb{Z}_2)$ one has $\|\alpha\|^{\mathrm{im}}_{\mathbb{Z}_2}=0$ while $\|\alpha\|^{\mathrm{em}}_{\mathbb{Z}_2}=k-2$, the latter being the minimal non-orientable genus computed by Bredon and Wood \cite{BW}. Hence the gap between the two norms can be arbitrarily large.

When discussing $\|\alpha\|^{\mathrm{im}}$ and $\|\alpha\|^{\mathrm{im}}_{\mathbb{Z}_2}$, we use immersed surfaces $F \looparrowright M$ in the smooth category. By \cite[Lemma 1.2]{HassHughes}, a continuous map $F \to M$ can always be homotoped to an immersion, so we may just as well work in the continuous category. Moreover, in most cases the problem of finding an immersion can be reduced to that of finding a homomorphism between fundamental groups. Indeed, when $M$ is orientable and irreducible we have $\pi_2(M)=0$, and since $F$ is two-dimensional this allows us to produce a continuous map $F \to M$ from any homomorphism $\pi_1(F) \to \pi_1(M)$. Hence, in order to calculate the immersed norm of a class in $H_2(M;\mathbb{Z}_2)$, we need only find a suitable homomorphism of fundamental groups. In Section \ref{sec-prep} we reduce the problem further, to the study of the order 2 elements of $H_1(M;\mathbb{Z})$.

Among the eight Thurston geometries, exactly seven give non-hyperbolic closed orientable geometric 3-manifolds. Six of these classes are Seifert fibered and the remaining one is Sol. The $\mathbb{Z}_2$-Thurston norm $\|\cdot\|^{\mathrm{em}}_{\mathbb{Z}_2}$ has been determined for orientable Seifert fibered manifolds with orientable base \cite{Du} and for Sol manifolds \cite{DuWang}. The immersed norm $\|\cdot\|^{\mathrm{im}}_{\mathbb{Z}_2}$, by contrast, has not been determined for any of these manifolds. This paper is devoted to the calculation of $\|\cdot\|^{\mathrm{im}}_{\mathbb{Z}_2}$ for Seifert fibered manifolds with orientable or non-orientable base, and for Sol manifolds.

As for the immersed norm, the only results we are aware of are the ones in \cite{BW} and \cite{HughesKim}. The closed 3-manifolds Bredon and Wood considered are lens spaces. The closed 3-manifolds Hughes and Kim considered are obtained by Dehn surgery on torus knots. Hughes and Kim also considered the analogous relative problem for knot complements, where the relevant invariant is the crosscap number of a knot. We refer to Section \ref{sec-imnorm} (Remark \ref{rem-crosscap}) for a comparison with our results.

\subsection*{Main results}

The main results of our paper are as follows.

\begin{Theorem}\label{thm-main-orientable}
Let $\Sigma_g$ be a closed orientable surface of genus $g$, and let $M=(\Sigma_g,(p_1,q_1),\cdots,(p_n,q_n))$ be a closed orientable Seifert fibered manifold with $g \geqslant 1$ or $n \geqslant 3$, where $(p_i,q_i)$ are the Seifert invariants. Then the structure of $H_2(M;\mathbb{Z}_2)$ and the values of $\|\cdot\|^{\mathrm{im}}_{\mathbb{Z}_2}$ on the elements of $H_2(M;\mathbb{Z}_2)$ are as follows:
\begin{enumerate}
\item[(1)] if $p_1,\cdots,p_n$ are all odd, then
\begin{enumerate}
\item[(1.1)] when $q_1 + \cdots + q_n$ is odd, $H_2(M;\mathbb{Z}_2)$ is isomorphic to $\mathbb{Z}_2^{2g}$, and every element $\alpha$ of it satisfies $\|\alpha\|^{\mathrm{im}}_{\mathbb{Z}_2}=0$;

\item[(1.2)] when $q_1 + \cdots + q_n$ is even, $H_2(M;\mathbb{Z}_2)$ contains a subgroup $L \cong \mathbb{Z}_2^{2g}$ (see Section \ref{sec-orientable}) such that every $\alpha \in L$ satisfies $\|\alpha\|^{\mathrm{im}}_{\mathbb{Z}_2}=0$ and every $\beta \in H_2(M;\mathbb{Z}_2)\setminus L$ satisfies $\|\beta\|^{\mathrm{im}}_{\mathbb{Z}_2} = 2g+n-2$;
\end{enumerate}

\item[(2)] if at least one of $p_1,\cdots,p_n$ is even, then every element $\alpha\in H_2(M;\mathbb{Z}_2)$ satisfies $\|\alpha\|^{\mathrm{im}}_{\mathbb{Z}_2}=0$.
\end{enumerate}
\end{Theorem}

\begin{Remark}\label{rem-exclude}
The hypothesis $g \geqslant 1$ or $n \geqslant 3$ in Theorem \ref{thm-main-orientable} serves to exclude $S^2 \times S^1$ and the lens spaces. For $S^2 \times S^1$ the result is obvious. For lens spaces the results are known. See the discussion following Proposition \ref{prop-h2-orientable}.
\end{Remark}

\begin{Theorem}\label{thm-main-nonorientable}
Let $N_g$ be the connected sum of $g$ copies of $\mathbb{R}P^2$, and let $M=(N_g,(p_1,q_1),\cdots,(p_n,q_n))$ be a closed orientable Seifert fibered manifold with $g+n \geqslant 3$, where $(p_i,q_i)$ are the Seifert invariants. Then the structure of $H_2(M;\mathbb{Z}_2)$ and the values of $\|\cdot\|^{\mathrm{im}}_{\mathbb{Z}_2}$ on the elements of $H_2(M;\mathbb{Z}_2)$ are as follows:
\begin{enumerate}
\item[(1)] if $p_1,\cdots,p_n$ are all odd, then
\begin{enumerate}
\item[(1.1)] when $q_1 + \cdots + q_n$ is odd, $H_2(M;\mathbb{Z}_2)$ is isomorphic to $\mathbb{Z}_2^{g}$, and every element $\alpha$ of it satisfies $\|\alpha\|^{\mathrm{im}}_{\mathbb{Z}_2}=0$;

\item[(1.2)] when $q_1 + \cdots + q_n$ is even, $H_2(M;\mathbb{Z}_2)$ contains a subgroup $L \cong \mathbb{Z}_2^{g}$ (see Section \ref{sec-nonorientable}) such that every $\alpha \in L$ satisfies $\|\alpha\|^{\mathrm{im}}_{\mathbb{Z}_2}=0$ and every $\beta \in H_2(M;\mathbb{Z}_2)\setminus L$ satisfies $\|\beta\|^{\mathrm{im}}_{\mathbb{Z}_2} = g+n-2$;
\end{enumerate}

\item[(2)] if at least one of $p_1,\cdots,p_n$ is even, then every element $\alpha\in H_2(M;\mathbb{Z}_2)$ satisfies $\|\alpha\|^{\mathrm{im}}_{\mathbb{Z}_2}=0$.
\end{enumerate}
\end{Theorem}

\begin{Remark}\label{rem-smallg}
If $g$ and $n$ do not satisfy $g+n \geqslant 3$ in Theorem \ref{thm-main-nonorientable}, then $M$ admits another Seifert fibration with orientable base (see \cite[Section 10.4]{Martelli}), hence falls into the cases covered by Theorem \ref{thm-main-orientable}.
\end{Remark}

Every Sol manifold is either a torus bundle over the circle or a torus semi-bundle, and both are obtained from an Anosov map of $T^2$.

\begin{Theorem}\label{thm-torusbundle}
Let $M$ be the mapping torus of an Anosov map of $T^2$. Then $H_2(M;\mathbb{Z}_2)$ is isomorphic to $\mathbb{Z}_2$, $\mathbb{Z}_2^2$, or $\mathbb{Z}_2^3$. In all cases, $H_2(M;\mathbb{Z}_2)$ contains a subgroup $L \cong \mathbb{Z}_2$ such that the non-trivial element $\alpha \in L$ satisfies $\|\alpha\|^{\mathrm{im}}_{\mathbb{Z}_2}=0$ and every $\beta \in H_2(M;\mathbb{Z}_2)\setminus L$ satisfies $\|\beta\|^{\mathrm{im}}_{\mathbb{Z}_2}=1$.
\end{Theorem}

\begin{Theorem}\label{thm-semibundle}
Let $M$ be the torus semi-bundle obtained from an Anosov map of $T^2$. Then $H_2(M;\mathbb{Z}_2)$ is isomorphic to $\mathbb{Z}_2^2$ or $\mathbb{Z}_2^3$. In both cases, $H_2(M;\mathbb{Z}_2)$ contains a subgroup $L\cong \mathbb{Z}_2^2$ such that every $\alpha \in L$ satisfies $\|\alpha\|^{\mathrm{im}}_{\mathbb{Z}_2}=0$ and every $\beta \in H_2(M;\mathbb{Z}_2)\setminus L$ satisfies $\|\beta\|^{\mathrm{im}}_{\mathbb{Z}_2}=1$.
\end{Theorem}

When the immersed norm vanishes and the simplest representing surface is connected and non-orientable, the definition does not distinguish between a Klein bottle and a projective plane. In Section \ref{sec-imnorm} we determine all the cases in which the simplest immersed representing surfaces are projective planes.

Bredon and Wood in \cite{BW} already found that the simplest embedded non-orientable representing surface of a class in $H_2(M;\mathbb{Z}_2)$ need not be $\pi_1$-injective. As we shall see (Section \ref{sec-imnorm}), even after passing from embedded to immersed surfaces, the simplest immersed non-orientable representing surface may still fail to be $\pi_1$-injective.

\subsection*{Organization of the Paper}

Section \ref{sec-prep} proves several facts that will be used in the later sections. Section \ref{sec-orientable} treats Seifert fibered manifolds with orientable base surface. Section \ref{sec-nonorientable} treats Seifert fibered manifolds with non-orientable base surface. Section \ref{sec-torusbundle} treats torus bundles. Section \ref{sec-semibundle} treats torus semi-bundles. Section \ref{sec-imnorm} determines all the cases in which the simplest immersed representing surfaces are projective planes.

\subsection*{Notation and Conventions}

Throughout the paper, $T^2$ denotes the torus and $K$ the Klein bottle $\mathbb{R}P^2 \# \mathbb{R}P^2$. By a punctured M\"obius band, we mean a M\"obius band with one (open) disc removed, hence such a surface is a compact non-orientable $2$-manifold with two boundary components.

\subsection*{Acknowledgements}

We thank Xiaoyu Xu for pointing out that we had missed the case of Seifert fibered manifolds with non-orientable base in the first draft of the paper.

\section{Preliminaries}\label{sec-prep}

\subsection{Detecting a 2-dimensional mod 2 homology class by an immersed projective plane}\label{subsec-detect}

By the universal coefficient theorem, $H_2(M;\mathbb{Z}_2)$ fits into the exact sequence
$$
0 \to H_2(M;\mathbb{Z})\otimes \mathbb{Z}_2 \to H_2(M;\mathbb{Z}_2) \to \operatorname{Tor}(H_1(M;\mathbb{Z}),\mathbb{Z}_2) \to 0,
$$
where $\operatorname{Tor}(G,\mathbb{Z}_2)$ denotes the subgroup consisting of all elements of order $2$ in $G$. In other words, we have $H_2(M;\mathbb{Z}_2) \cong (H_2(M;\mathbb{Z})\otimes \mathbb{Z}_2) \oplus \operatorname{Tor}(H_1(M;\mathbb{Z}),\mathbb{Z}_2)$. The elements of the first summand can be represented by orientable surfaces, whereas every element outside $H_2(M;\mathbb{Z})\otimes \mathbb{Z}_2$ can only be represented by non-orientable surfaces.

\begin{Proposition}\label{prop-pi2-detection}
Let $M$ be an orientable irreducible 3-manifold. There is an immersion $i : \mathbb{R}P^2 \looparrowright M$ with $[i(\mathbb{R}P^2)] \notin H_2(M;\mathbb{Z})\otimes \mathbb{Z}_2$ if and only if $\pi_1(M)$ contains an element of order $2$ which, after abelianization, becomes an element of order $2$ in $H_1(M;\mathbb{Z})$.
\end{Proposition}

\begin{proof}
For a map $f : N \to M$, the Bockstein maps $\partial_* : H_2(-;\mathbb{Z}_2) \to H_1(-;\mathbb{Z})$ for $N$ and for $M$ are natural, that is, $f_* \circ \partial_* = \partial_* \circ f_*$. For the projective plane, $\partial_* : H_2(\mathbb{R}P^2;\mathbb{Z}_2) \to H_1(\mathbb{R}P^2;\mathbb{Z})$ is an isomorphism onto $H_1(\mathbb{R}P^2;\mathbb{Z}) \cong \mathbb{Z}_2$, and it sends the nonzero element of $H_2(\mathbb{R}P^2;\mathbb{Z}_2)$ to the element of order $2$.

Suppose first that $\pi_1(M)$ contains an element $x$ of order $2$ whose image in $H_1(M;\mathbb{Z})$ has order $2$. Sending the generator $g$ of $\pi_1(\mathbb{R}P^2) = \mathbb{Z}_2$ to $x$ defines a homomorphism $\pi_1(\mathbb{R}P^2) \to \pi_1(M)$. Since $M$ is irreducible, $\pi_2(M) = 0$, so this homomorphism is realized by a map $f : \mathbb{R}P^2 \to M$, and by \cite[Lemma 1.2]{HassHughes} we may take $f$ to be an immersion. By naturality, $\partial_*[f(\mathbb{R}P^2)] = f_*(\partial_*[\mathbb{R}P^2]) = x$ in $H_1(M;\mathbb{Z})$. This image is nonzero, so $[f(\mathbb{R}P^2)]$ lies outside $H_2(M;\mathbb{Z})\otimes \mathbb{Z}_2$.

Conversely, let $i : \mathbb{R}P^2 \looparrowright M$ be an immersion whose class $\alpha = [i(\mathbb{R}P^2)]$ lies outside $H_2(M;\mathbb{Z})\otimes \mathbb{Z}_2$. Then $\partial_*\alpha \neq 0$. Let $x \in \pi_1(M)$ be the image of the generator $g$ of $\pi_1(\mathbb{R}P^2)$ under $i_\pi$. Since $g^2 = 1$, we get $x^2 = i_\pi(g^2) = 1$, so $x$ has order at most $2$. By naturality, the image of $x$ in $H_1(M;\mathbb{Z})$ equals $i_*(\partial_*[\mathbb{R}P^2]) = \partial_*\alpha$, which is nonzero. Hence $x$ has order exactly $2$.

The assumption of irreducibility is used only to have $\pi_2(M)=0$ and hence to realize a homomorphism $\pi_1(\mathbb{R}P^2)\to\pi_1(M)$ by a map. It is therefore unnecessary when $\pi_1(M)=\mathbb{Z}$, in particular for $M\cong S^2\times S^1$, which is the only reducible manifold to which the result is applied. There $\pi_1(M)$ has no element of order $2$ and $H_2(M;\mathbb{Z}_2)=H_2(M;\mathbb{Z})\otimes\mathbb{Z}_2$, so neither side of the equivalence holds.
\end{proof}

\begin{Remark}\label{rem-other-nonorientable}
For a non-orientable surface $N$ other than $\mathbb{R}P^2$ the analogue of Proposition \ref{prop-pi2-detection} fails. Let $\alpha \in H_2(M;\mathbb{Z}_2) \setminus H_2(M;\mathbb{Z})\otimes \mathbb{Z}_2$. A homomorphism $\pi_1(N) \to \pi_1(M)$ that sends the unique element of order $2$ in $H_1(N;\mathbb{Z})$ to $\partial_*\alpha$ is realized by a surface representing $\alpha$ plus some element of $H_2(M;\mathbb{Z})\otimes \mathbb{Z}_2$, not necessarily by one representing $\alpha$ itself, because $\partial_*$ factors out $H_2(M;\mathbb{Z})\otimes \mathbb{Z}_2$. For Seifert fibered manifolds such homomorphisms can be found by Gr\"obner--Shirshov bases \cite{LiaoZhao}, but in general their existence is undecidable.
\end{Remark}

\subsection{Immersed M\"obius bands in a solid torus}\label{subsec-mobius-solid}

For the solid torus $D^2 \times S^1$, a curve on its boundary torus that runs $n$ times in the $\partial D^2$ direction and $m$ times in the $S^1$ direction is called an $(m,n)$-curve. Here the $S^1$ direction is the longitude of the solid torus and the $\partial D^2$ direction is its meridian.

\begin{Proposition}\label{prop-mobius-solid}
Let $N$ be a solid torus and $S$ a M\"obius band. For coprime $2k$ and $q$, there is an immersion $f: (S,\partial S) \to (N,\partial N)$ such that $f(\partial S)$ is the $(2k,q)$-curve on $\partial N$.
\end{Proposition}

\begin{figure}[htbp]
\begin{center}
\includegraphics[width=0.55\textwidth]{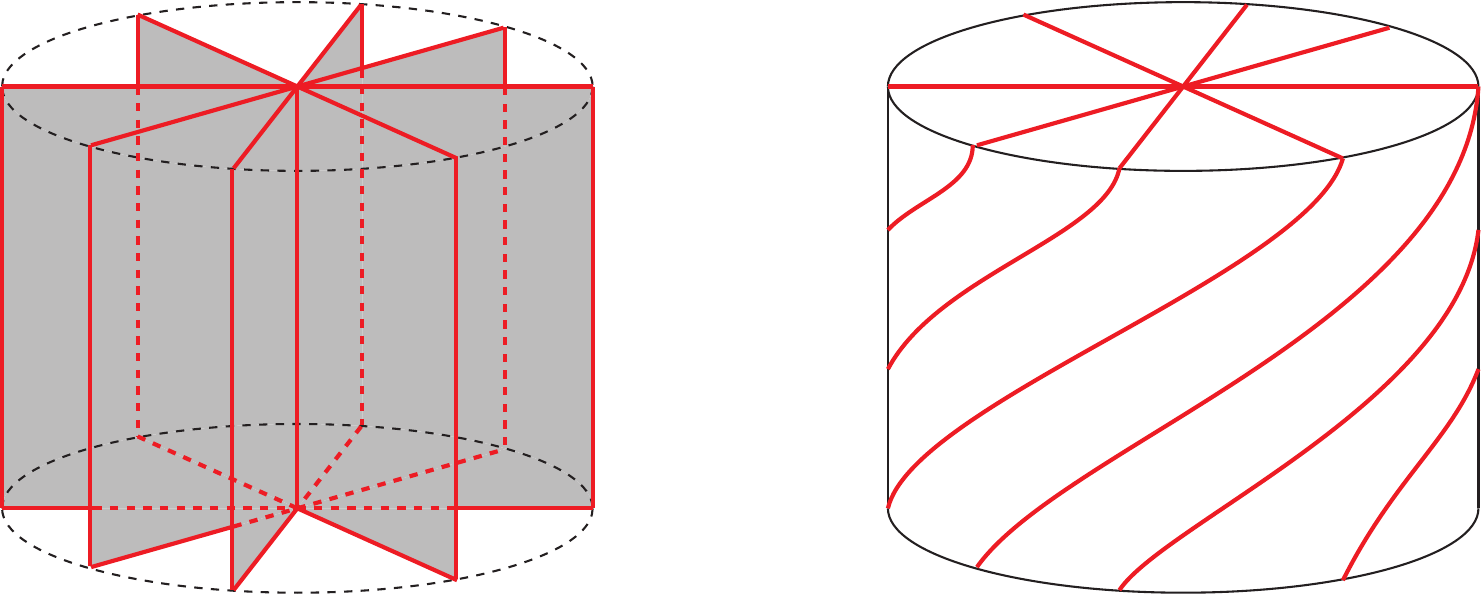}
\caption{The case $k=4$. Rotate the upper base of the cylinder counterclockwise through $3\pi/4$ and glue it to the lower base. The result is an immersed M\"obius band in the solid torus whose boundary is the $(8,3)$-curve on $\partial (D^2 \times S^1)$.}
\label{fig-solid-torus}
\end{center}
\end{figure}

\begin{proof}
(Following the construction in \cite[Section 4.4]{Carter}.) Let $D^2$ be the unit disk, and let $X \subset D^2$ be the union of the $2k$ rays issuing from the origin at angles $\theta_i=i\pi/k$, $i=0,1,\cdots,2k-1$. Denote by $r_{\varphi}$ the counterclockwise rotation of $D^2$ through the angle $\varphi$. In $D^2 \times [0,1]$ take the subset
$$
Y=\{(x,t)\mid x \in r_{\pi q t/k}(X)\subset D^2,\; t \in [0,1]\}.
$$
After gluing $D^2 \times \{0\}$ to $D^2 \times \{1\}$, which yields the solid torus $D^2 \times S^1$, the set $Y$ becomes an immersed M\"obius band in $D^2 \times S^1$ whose boundary is the curve on $\partial(D^2 \times S^1)$ that runs $2k$ times in the $S^1$ direction and $q$ times in the $\partial D^2$ direction (hereafter the $(2k,q)$-curve). See Figure \ref{fig-solid-torus}.
\end{proof}

\subsection{Immersed punctured M\"obius bands in $T^2 \times [0,1]$}\label{subsec-punctured}

Recall that a punctured M\"obius band means a M\"obius band with one (open) disc removed. For $T^2 \times I$, given two curves $c_0 \subset T^2 \times 0$ and $c_1 \subset T^2 \times 1$, we look for the minimal genus of a non-orientable surface immersed in $T^2 \times I$ with two boundary curves $c_0$ and $c_1$.

\begin{Proposition}\label{prop-punctured-band}
For $T^2 \times I$, let $c_0$ be a curve on $T^2 \times 0$ and $c_1$ a curve on $T^2 \times 1$, and let $p(c_0),p(c_1)$ denote the projections of $c_0,c_1$ to the torus $T^2$. If a punctured M\"obius band is immersed in $T^2 \times I$ with boundary $c_0 \cup c_1$, then the mod 2 intersection number of $p(c_0)$ and $p(c_1)$ is $0$. Conversely, if that intersection number is $0$ and $p(c_1)\neq \pm p(c_0)$, then such a band exists.
\end{Proposition}

\begin{proof}
Necessity follows from the fact that $p(c_0)$ and $p(c_1)$ lie in the same mod 2 homology class of $T^2$. Sufficiency follows from the construction below.

We may assume that $p(c_0)$ is the $(0,1)$-curve on the torus and $p(c_1)$ is the $(2k,q)$-curve on the torus, with $k\geqslant 1$, which is possible because $p(c_1)\neq \pm p(c_0)$. Let $D^2 \setminus D^2_{\varepsilon}$ be the standard disk with a concentric open disk of radius $\varepsilon$ removed. Then $T^2 \times I$ is homeomorphic to $(D^2 \setminus D^2_{\varepsilon})\times S^1$, where $T^2 \times 0$ corresponds to the inner boundary $\partial D^2_{\varepsilon} \times S^1$ and $T^2 \times 1$ to the outer boundary $\partial D^2 \times S^1$. We construct directly a punctured M\"obius band immersed in $(D^2 \setminus D^2_{\varepsilon})\times S^1$ whose boundary consists of the $(2k,q)$-curve on $\partial D^2 \times S^1$ and the $(0,1)$-curve on $\partial D^2_{\varepsilon} \times S^1$.

\begin{figure}[htbp]
\begin{center}
\includegraphics[width=0.9\textwidth]{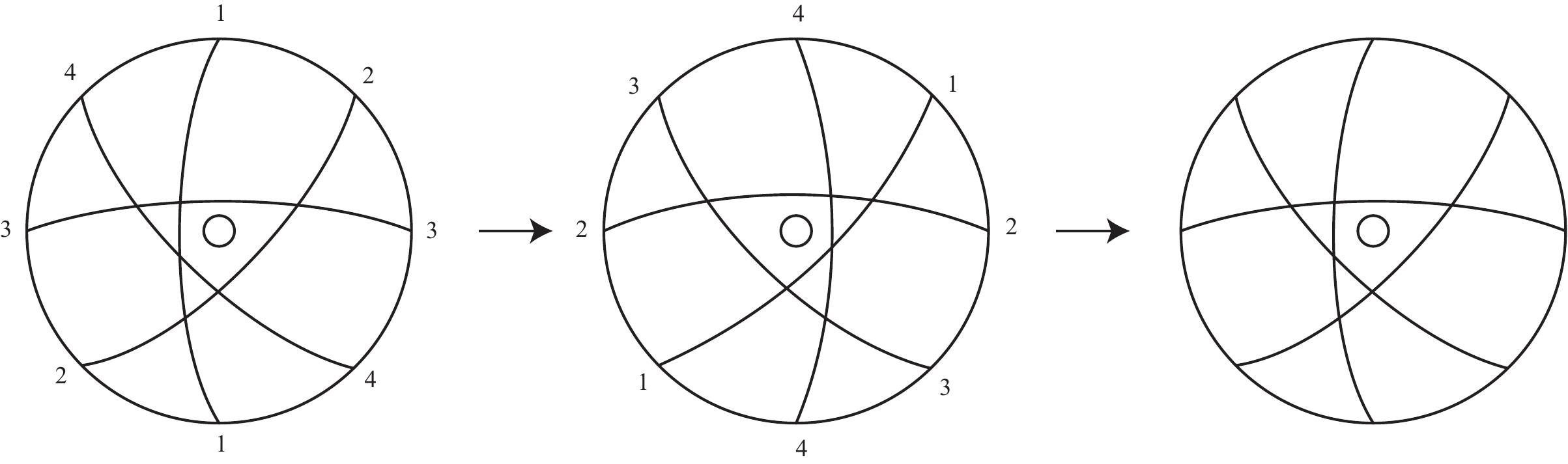}
\caption{The arcs $\alpha_1,\cdots,\alpha_k$ in the annulus $D^2 \setminus D^2_{\varepsilon}$, for $k=4$ and $q=3$.}
\label{fig-arcs}
\end{center}
\end{figure}

Choose an arc $\alpha_1 \subset D^2$ joining a pair of antipodal points on $\partial D^2$ and avoiding $D^2_{\varepsilon}$. Next, construct $k-1$ further arcs in $D^2 \setminus D^2_{\varepsilon}$: for $i=1,\cdots,k-1$, let $\alpha_{1+i}$ be the image of $\alpha_1$ under the counterclockwise rotation through the angle $iq\pi/k$. This yields $\alpha_2,\cdots,\alpha_k$. The leftmost part of Figure \ref{fig-arcs} shows the case $k=4$, $q=3$. Denote the union of all these arcs by $X'$.

Denote by $r_{\varphi}$ the counterclockwise rotation of $D^2 \setminus D^2_{\varepsilon}$ through the angle $\varphi$. By construction we have $r_{q\pi/k}(\alpha_i)=\alpha_{i+1}$ ($i=1,2,\cdots,k-1$), and $r_{q\pi/k}(\alpha_k)$ is the image of $\alpha_1$ under the central symmetry about the origin.

In $(D^2 \setminus D^2_{\varepsilon}) \times [0,\frac 12]$, take the subset
$$
Y'=\{(x,t)\mid x \in r_{{2q\pi t}/k}(X')\subset D^2 \setminus D^2_{\varepsilon},\;t \in [0,\textstyle\frac 12]\}.
$$
Equivalently, $Y'$ is obtained from $X' \times [0,\frac 12]$ by twisting the upper face of $(D^2 \setminus D^2_{\varepsilon})\times [0,\frac 12]$ counterclockwise through the angle $q\pi/k$. The left part of Figure \ref{fig-punctured} is $X' \times [0,\frac 12]$ in $(D^2 \setminus D^2_{\varepsilon})\times [0,\frac 12]$, and the middle part of Figure \ref{fig-punctured} is $Y'$ in $(D^2 \setminus D^2_{\varepsilon})\times [0,\frac 12]$.

\begin{figure}[htbp]
\begin{center}
\includegraphics[width=0.9\textwidth]{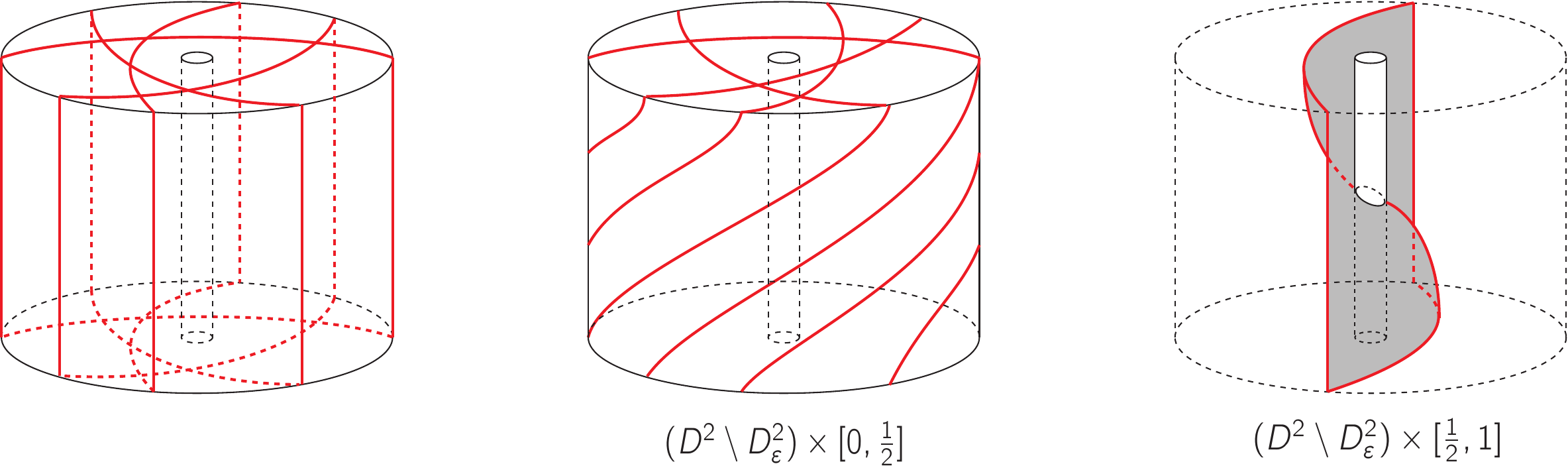}
\caption{From left to right: $X' \times [0,\frac 12]$ in $(D^2 \setminus D^2_{\varepsilon})\times [0,\frac 12]$, the twisted surface $Y'$, and the surface $Y''$ in $(D^2 \setminus D^2_{\varepsilon}) \times [\frac 12,1]$.}
\label{fig-punctured}
\end{center}
\end{figure}

The restrictions of $Y'$ to $(D^2 \setminus D^2_{\varepsilon}) \times 0$ and to $(D^2 \setminus D^2_{\varepsilon}) \times \frac 12$ differ only on the arc $\alpha_1$. Next we construct an immersed surface $Y''$ in $(D^2 \setminus D^2_{\varepsilon}) \times [\frac 12,1]$. Take a homotopy $F_t: \alpha_1 \times [0,1] \to D^2$ from $r_{\pi}\alpha_1$ to $\alpha_1$, satisfying $F_0(\alpha_1)=r_{\pi}\alpha_1$ and $F_1(\alpha_1)=\alpha_1$, and such that $F_t$ sweeps over $D^2_{\varepsilon}$ exactly once. In $D^2 \times [\frac 12,1]$ take $\{(x,\frac 12 + \frac 12t) \mid x \in F_t(\alpha_1), t \in [0,1]\}$, and then remove $D^2_{\varepsilon} \times [\frac 12,1]$. See the rightmost part of Figure \ref{fig-punctured}. Taking the union with $(\alpha_2 \cup \cdots \cup \alpha_k) \times [\frac 12,1]$, we obtain an immersed surface $Y''$ in $(D^2 \setminus D^2_{\varepsilon}) \times [\frac 12,1]$. Now $Y''\cap ((D^2 \setminus D^2_{\varepsilon}) \times \frac 12)$ coincides with $Y' \cap ((D^2 \setminus D^2_{\varepsilon}) \times \frac 12)$, and $Y'' \cap ((D^2 \setminus D^2_{\varepsilon}) \times 1)$ coincides with $Y' \cap ((D^2 \setminus D^2_{\varepsilon}) \times 0)$. Hence, gluing $(D^2 \setminus D^2_{\varepsilon}) \times [0,\frac 12]$ to $(D^2 \setminus D^2_{\varepsilon}) \times [\frac 12,1]$ along $(D^2 \setminus D^2_{\varepsilon}) \times \frac 12$, and then identifying $(D^2 \setminus D^2_{\varepsilon}) \times 0$ with $(D^2 \setminus D^2_{\varepsilon}) \times 1$, we obtain a punctured M\"obius band immersed in $(D^2 \setminus D^2_{\varepsilon}) \times S^1$. The boundary of such an immersed punctured M\"obius band consists of a $(2k,q)$-curve on $\partial D^2 \times S^1$ and a $(0,1)$-curve on $\partial D^2_{\varepsilon} \times S^1$.
\end{proof}

\subsection{Classification of the nonzero elements of $H_1(T^2;\mathbb{Z}_2)$ by mod 2 \mbox{intersection} number}\label{subsec-three-classes}

\begin{Lemma}\label{lem-three-classes}
For classes $c,c'\in H_1(T^2;\mathbb{Z}_2)$, let $i_2(c,c')$ denote their mod 2 intersection number. On the set of nonzero elements of $H_1(T^2;\mathbb{Z}_2)$, the relation $i_2(c,c')=0$ is an equivalence relation with exactly $3$ equivalence classes: the $(1,0)$-class, the $(0,1)$-class, and the $(1,1)$-class.
\end{Lemma}

\begin{proof}
It suffices to compute the determinants of the corresponding $2 \times 2$ matrices of these curves. When two curves lie in different classes, the mod 2 residue class of the value of the determinant is 1. When two curves lie in the same class, the mod 2 residue class of the value of the determinant is $0$.
\end{proof}

\section{Seifert Fibered Manifolds with Orientable Base}\label{sec-orientable}

Let $M$ be a closed Seifert fibered manifold whose base surface is orientable and whose total space is orientable. We write $M=(\Sigma_g,(p_1,q_1),\cdots,(p_n,q_n))$, where $\Sigma_g$ is the closed orientable surface of genus $g$. It is obtained by removing $n$ small disks $D(P_1),\cdots,D(P_n)$ around the singular points $P_1,\cdots,P_n$ of the base orbifold, which gives $B_0=\Sigma_g - nD^2$, and then gluing $n$ solid tori $N_1,\cdots,N_n$ onto the boundary of $M_0=B_0 \times S^1$. The gluing map of $N_i$ identifies the meridian of $N_i$ with the curve on $\partial B_0\times S^1$ that runs $p_i$ times in the $\partial D(P_i)$ direction and $q_i$ times in the $S^1$ direction. Here $p_i$ and $q_i$ are coprime with $p_i>0$, so $p_i$ is the multiplicity of the $i$-th singular fiber and the base orbifold has a cone point of angle $2\pi/p_i$ at $P_i$. Thus $p_i\geqslant 2$, hence $q_i\neq 0$, whenever $n\geqslant 1$. The gluing maps determine $M$, and all computations below use the data $(p_i,q_i)$ as given. A fibration without singular fibers is recorded in the same way by the single pair $(p,q)=(1,e)$, where $e$ is the characteristic class of the underlying circle bundle. For the product $\Sigma_g\times S^1$ this is the pair $(1,0)$.

For $M=(\Sigma_g,(p_1,q_1),\cdots,(p_n,q_n))$, a straightforward computation gives its fundamental group and its homology groups as follows:
$$
\begin{array}{l}
\pi_1(M) = \langle a_1,b_1,\cdots,a_g,b_g,h_1,\cdots,h_n,v \mid \\
\qquad a_iv=va_i, \quad b_iv=vb_i \quad (i = 1,\cdots,g), \\
\qquad h_jv = vh_j, \quad h_j^{p_j}v^{q_j} = 1 \quad (j=1,\cdots,n), \\
\qquad h_1\cdots h_n[a_1,b_1]\cdots [a_g,b_g] = 1 \rangle
\end{array}
$$
$$
\begin{array}{l}
H_1(M;\mathbb{Z}) = \langle [a_1],[b_1],\cdots,[a_g],[b_g],[h_1],\cdots,[h_n],[v] \;\; \mid \\
\quad p_j[h_j] + q_j[v]=0\ (j=1,\cdots,n), \quad [h_1] + \cdots + [h_n]=0 \rangle
\end{array}
$$

Below, Proposition \ref{prop-h2-orientable} computes $H_2(M;\mathbb{Z}_2)$, Proposition \ref{prop-upper-bound-orientable} gives upper bounds for the immersed norm, and Proposition \ref{prop-lower-bound-orientable} shows that these bounds are attained.

\begin{Proposition}\label{prop-h2-orientable}
Let $M = (\Sigma_g,(p_1,q_1),\cdots,(p_n,q_n))$ be a Seifert fibered manifold. Then:
\begin{enumerate}
\item[(1)] if $p_1,\cdots,p_n$ are all odd, then
\begin{enumerate}
\item[(1.1)] when $q_1+\cdots+q_n$ is odd, we have $H_2(M;\mathbb{Z}_2) \cong \mathbb{Z}_2^{2g}$;

\item[(1.2)] when $q_1+\cdots+q_n$ is even, we have $H_2(M;\mathbb{Z}_2) \cong \mathbb{Z}_2^{2g+1}$;
\end{enumerate}

\item[(2)] if among $p_1,\cdots,p_n$ there are $k \geqslant 1$ even numbers, then we have $H_2(M;\mathbb{Z}_2) \cong \mathbb{Z}_2^{2g+k-1}$.
\end{enumerate}
\end{Proposition}

\begin{proof}
Since $M$ is closed and orientable, Poincar\'e duality identifies $H_2(M;\mathbb{Z}_2)$ with $H_1(M;\mathbb{Z}_2)$. The latter is obtained from the presentation of $H_1(M;\mathbb{Z})$ above by reducing all coefficients modulo 2, an even coefficient becoming $0$ and an odd one becoming $1$. The classes $[a_i]$ and $[b_i]$ occur in no relation and contribute $2g$ free generators. It remains to see what the reduction does to $[h_1],\cdots,[h_n]$ and $[v]$.

If $p_1,\cdots,p_n$ are all odd, then the reduced relation $p_j[h_j]+q_j[v]=0$ reads $[h_j]=q_j[v]$, so every $[h_j]$ is determined by $[v]$. The relation $[h_1]+\cdots+[h_n]=0$ then reads $(q_1+\cdots+q_n)[v]=0$. Hence $[v]$, and with it all the $[h_j]$, survives exactly when $q_1+\cdots+q_n$ is even, and the dimension is $2g$ or $2g+1$.

If $k\geqslant 1$ of the $p_j$ are even, take one such $j$. Since $\gcd(p_j,q_j)=1$, the number $q_j$ is odd, so the reduced relation $p_j[h_j]+q_j[v]=0$ reads $[v]=0$. With $[v]=0$, each relation with odd $p_j$ forces $[h_j]=q_j[v]=0$, so the only classes that survive are the $k$ classes $[h_j]$ with $p_j$ even. These $k$ classes are still subject to the sum relation $[h_1]+\cdots+[h_n]=0$, and every other term in it is now zero, so the relation says precisely that the sum of these $k$ classes is zero. No other relation involves them, so one relation among $k$ classes leaves $k-1$ of them independent, and the dimension is $2g+k-1$.
\end{proof}

We denote by $[v]$ the class of a regular fiber in $H_1(M;\mathbb{Z}_2)$, and by $i_2(\alpha,\gamma)$ the mod 2 intersection number of $\alpha\in H_2(M;\mathbb{Z}_2)$ with $\gamma\in H_1(M;\mathbb{Z}_2)$. We define $L\subset H_2(M;\mathbb{Z}_2)$ to be the subgroup generated by the classes of the $2g$ vertical tori $a_i\times S^1$ and $b_i\times S^1$ ($i=1,\cdots,g$), where $a_1,b_1,\cdots,a_g,b_g$ are the curves on $\Sigma_g$ appearing in the presentation of $\pi_1(M)$ above. Since $i_2([a_i\times S^1],[b_j])=\delta_{ij}$, these $2g$ classes are independent, so $L\cong\mathbb{Z}_2^{2g}$, and after small isotopies a sum of classes in $L$ is represented by a disjoint union of such tori, so every $\alpha\in L$ satisfies $\|\alpha\|^{\mathrm{im}}_{\mathbb{Z}_2}=0$. The pairing with $[v]$ gives a homomorphism $\psi(\alpha)=i_2(\alpha,[v])$ from $H_2(M;\mathbb{Z}_2)$ to $\mathbb{Z}_2$, and each torus generating $L$ is disjoint from the fibers over the points of $\Sigma_g$ outside the corresponding curve, so $L\subset\ker\psi$.

When $g=0$ and $n=0$, we have $M\cong S^2\times S^1$. Here $H_2(M;\mathbb{Z}_2)\cong\mathbb{Z}_2$, and the nonzero element is represented by the sphere factor $S^2\times\{\mathrm{pt}\}$, so its immersed norm is $0$.

When $g=0$ and $1\leqslant n\leqslant 2$, $M$ is a lens space $L(p,q)$, and $H_2(L(p,q);\mathbb{Z}_2)\cong H_1(L(p,q);\mathbb{Z})\otimes\mathbb{Z}_2$ is nonzero only when $p$ is even. In that case $M=L(2k,q)$, and the construction of immersed M\"obius bands in solid tori from Section \ref{subsec-mobius-solid} gives an immersed projective plane in $L(2k,q)$ representing the nonzero element, so its immersed norm is $0$.

When $n=0$ and $g\geqslant 1$, $M$ is the product $\Sigma_g\times S^1$, and with the empty list recorded as $(1,0)$ Proposition \ref{prop-h2-orientable} gives $H_2(M;\mathbb{Z}_2)\cong\mathbb{Z}_2^{2g+1}$. The horizontal surface $\Sigma_g\times\{\mathrm{pt}\}$ has Euler characteristic $2-2g$ and meets every regular fiber once, so $i_2([\Sigma_g\times\{\mathrm{pt}\}],[v])=1$ and $\psi:H_2(M;\mathbb{Z}_2) \to \mathbb{Z}_2$ is nonzero. As $L$ has index $2$ in $H_2(M;\mathbb{Z}_2)$, this gives $\ker\psi=L$, so the class lies outside $L$. Every $\beta\in H_2(M;\mathbb{Z}_2)\setminus L$ is of the form $[\Sigma_g\times\{\mathrm{pt}\}]+\alpha$ with $\alpha\in L$, and is therefore represented by a surface of complexity at most $2g-2$. Conversely, let $F$ represent a class $\beta\in H_2(M;\mathbb{Z}_2)\setminus L$. Since $i_2(\beta,[v])=1$, some component of $F$ meets each regular fiber an odd number of times, hence is mapped onto $\Sigma_g$ with nonzero mod 2 degree. By Lemma \ref{lem-degree-bound}, proved below, a surface that maps onto $\Sigma_g$ with nonzero mod 2 degree has Euler characteristic at most $2-2g$. Hence this component contributes at least $2g-2$ to the complexity of $F$. Thus $\|\beta\|^{\mathrm{im}}_{\mathbb{Z}_2}=2g-2=2g+n-2$.

We may therefore assume that $M=(\Sigma_g,(p_1,q_1),\cdots,(p_n,q_n))$ satisfies $g\geqslant 1$ or $n\geqslant 3$, and that $n\geqslant 1$.

\begin{Proposition}\label{prop-upper-bound-orientable}
Let $M = (\Sigma_g,(p_1,q_1),\cdots,(p_n,q_n))$ be a Seifert fibered manifold with $g \geqslant 1$ or $n \geqslant 3$, and $n\geqslant 1$. Then:
\begin{enumerate}
\item[(1)] if $p_1,\cdots,p_n$ are all odd, then
\begin{enumerate}
\item[(1.1)] when $q_1+\cdots+q_n$ is odd, every element $\alpha$ of $H_2(M;\mathbb{Z}_2) \cong \mathbb{Z}_2^{2g}$ satisfies $\|\alpha\|^{\mathrm{im}}_{\mathbb{Z}_2}=0$;

\item[(1.2)] when $q_1+\cdots+q_n$ is even, every element $\alpha$ of $L$ satisfies $\|\alpha\|^{\mathrm{im}}_{\mathbb{Z}_2}=0$, and every element $\beta$ of $H_2(M;\mathbb{Z}_2)\setminus L$ satisfies $\|\beta\|^{\mathrm{im}}_{\mathbb{Z}_2} \leqslant 2g+n-2$;
\end{enumerate}

\item[(2)] if among $p_1,\cdots,p_n$ there are $k \geqslant 1$ even numbers, then every element $\alpha\in H_2(M;\mathbb{Z}_2)$ satisfies $\|\alpha\|^{\mathrm{im}}_{\mathbb{Z}_2}=0$.
\end{enumerate}
\end{Proposition}

\begin{proof}
We construct explicit immersed surfaces representing these mod 2 homology classes.

(1.1) When $p_1,\cdots,p_n$ are all odd and $q_1+\cdots+q_n$ is odd, we have $H_2(M;\mathbb{Z}_2)=L$, so there is nothing to prove.

(1.2) When $p_1,\cdots,p_n$ are all odd and $q_1+\cdots+q_n$ is even, we may assume that every $q_i$ is even. Indeed, replacing odd $q_i$ and odd $q_j$ by even $q_i+p_i$ and even $q_j-p_j$ only amounts to another choice of the partial cross-sections of the fibration and does not change $M$. As $q_1+\cdots+q_n$ is even, the odd $q_i$ come in pairs, so applying this to each pair makes every $q_i$ even.

The Seifert manifold $M$ is regarded as obtained from $(\Sigma_g - nD^2) \times S^1$ by gluing on $n$ solid tori. The surface $\Sigma_g-nD^2$ embeds naturally in $M$. The gluing map sends the meridian of $N_i$ to $p_i$ times the base direction plus $q_i$ times the fiber direction, with determinant $1$. Inverting it, one finds that each boundary circle of $\Sigma_g-nD^2$ wraps $-q_i$ times in the longitude direction of $\partial N_i$ (up to sign) and runs a number of times coprime to $q_i$ in the meridian direction. Since $q_i$ is even, Proposition \ref{prop-mobius-solid} applies, and each of these boundary circles bounds an immersed M\"obius band inside the corresponding glued solid torus. Gluing these $n$ M\"obius bands to $\Sigma_g-nD^2$, we get an immersed non-orientable surface $F$ of genus $2g+n$, whose Euler characteristic is $2-2g-n$.

Geometrically, the composite of the immersion $f: F \looparrowright M$ and the projection $\varphi: M \to \Sigma_g$ can be described as follows: it pinches each of the $n$ M\"obius bands to the corresponding singular point of the base orbifold. Since $i_2([f(F)],[v])=1$, $\psi:H_2(M;\mathbb{Z}_2) \to \mathbb{Z}_2$ is nonzero. As $L\subset\ker\psi$ has index $2$ in $H_2(M;\mathbb{Z}_2)$, $\ker\psi=L$. Now the class $\beta'=[f(F)]$ satisfies $\|\beta'\|^{\mathrm{im}}_{\mathbb{Z}_2} \leqslant 2g+n-2$. Every $\beta\in H_2(M;\mathbb{Z}_2)\setminus L$ has the form $\beta'+\alpha$ for some $\alpha\in L$. Since the norm vanishes on $L$, we have $\|\beta\|^{\mathrm{im}}_{\mathbb{Z}_2} \leqslant \|\beta'\|^{\mathrm{im}}_{\mathbb{Z}_2} \leqslant 2g+n-2$.

(2) When among $p_1,\cdots,p_n$ there are $k>0$ even ones, we may assume that $p_1,\cdots, p_k$ are even. Let $N_1,\cdots,N_k$ be the solid torus neighborhoods of these singular fibers. For any two distinct $i,j \in \{1,\cdots,k\}$, take a vertical annulus whose two boundary components lie on $\partial N_i$ and $\partial N_j$ respectively. By the inversion of the gluing map, the boundary circle of this annulus lying on $\partial N_i$ wraps $p_i$ times in the longitude direction of $\partial N_i$ and runs a number of times in the meridian direction that is coprime to $p_i$. The same description applies to the other boundary circle with $j$ in place of $i$. As $p_i$ and $p_j$ are even, Proposition \ref{prop-mobius-solid} applies: these two circles bound immersed M\"obius bands in $N_i$ and $N_j$ respectively, and together with the annulus they give a Klein bottle, whose Euler characteristic is $0$.

Counting all the Klein bottles obtained in this way together with all their disjoint unions, we obtain $\binom{k}{2} + \binom{k}{4} + \cdots + \binom{k}{2\lfloor k/2\rfloor} = 2^{k-1}-1$ surfaces in total. By computing their mod 2 intersection numbers with the simple closed curves that run horizontally, we see that they represent distinct homology classes, and these classes are exactly all the nonzero elements of the direct summand $\mathbb{Z}_2^{k-1}$ of $H_2(M;\mathbb{Z}_2)\cong L \oplus \mathbb{Z}_2^{k-1}$. Each of these elements $\beta$ satisfies $\|\beta\|^{\mathrm{im}}_{\mathbb{Z}_2}=0$.
\end{proof}

Thus it remains only to show that the inequality in Proposition \ref{prop-upper-bound-orientable} (1.2) is in fact an equality.

\begin{Lemma}\label{lem-degree-bound}
Let $\hat N$ and $\hat B$ be connected compact surfaces. Let $g:\hat N\to\hat B$ be a continuous map with $g(\partial\hat N)\subset\partial\hat B$. Suppose that $g$ maps the relative fundamental class of $\hat N$ to that of $\hat B$, that is, $g_*[\hat N,\partial\hat N]=[\hat B,\partial\hat B]$ in $H_2(-;\mathbb{Z}_2)$. Then $\chi(\hat N)\leqslant\chi(\hat B)$.
\end{Lemma}

\begin{proof}
For a connected compact surface $\hat F$ one has $b_1(\hat F;\mathbb{Z}_2)=1-\chi(\hat F)$. The cup product pairing $H^1(\hat F;\mathbb{Z}_2)\times H^1(\hat F,\partial\hat F;\mathbb{Z}_2)\to H^2(\hat F,\partial\hat F;\mathbb{Z}_2)\cong\mathbb{Z}_2$ is nondegenerate. Let $\alpha\in H^1(\hat B;\mathbb{Z}_2)$ be nonzero. Then there is $\beta\in H^1(\hat B,\partial\hat B;\mathbb{Z}_2)$ with $\langle\alpha\cup\beta,[\hat B,\partial\hat B]\rangle=1$. By naturality,
$$
\langle g^*\alpha\cup g^*\beta,[\hat N,\partial\hat N]\rangle=\langle\alpha\cup\beta,g_*[\hat N,\partial\hat N]\rangle=\langle\alpha\cup\beta,[\hat B,\partial\hat B]\rangle=1,
$$
Hence $g^*\alpha\neq 0$. Thus $g^*$ is injective. It follows that $b_1(\hat N;\mathbb{Z}_2)\geqslant b_1(\hat B;\mathbb{Z}_2)$, and therefore $\chi(\hat N)\leqslant\chi(\hat B)$.
\end{proof}

The mod 2 formulation needs no orientation, so it applies also to the non-orientable base of Section \ref{sec-nonorientable}.

\begin{Proposition}\label{prop-lower-bound-orientable}
Let $M = (\Sigma_g,(p_1,q_1),\cdots,(p_n,q_n))$ be a Seifert fibered manifold with $g\geqslant 1$ or $n\geqslant 3$, and $n\geqslant 1$, and suppose that all $p_i$ are odd, that no $q_i$ vanishes, and that $q_1+\cdots+q_n$ is even, so that $H_2(M;\mathbb{Z}_2) \cong \mathbb{Z}_2^{2g+1} = L \oplus \mathbb{Z}_2$. Then every element $\beta$ of $H_2(M;\mathbb{Z}_2)\setminus L$ satisfies $\|\beta\|^{\mathrm{im}}_{\mathbb{Z}_2} \geqslant 2g+n-2$.
\end{Proposition}

\begin{proof}
Let an immersed surface $f:N \looparrowright M$ represent $\beta \in H_2(M;\mathbb{Z}_2)\setminus L$, so that $\psi(\beta)=1$, the kernel of $\psi$ being $L$ in this case. If $N$ is not connected, apply the argument below to a component with nonzero $\psi$-value, which exists since $\psi$ is additive and $\psi(\beta)=1$. The other components contribute a nonnegative amount to the complexity. We may perturb $f$ so that it is transverse to a finite set of fibers which contains all the singular fibers and at least one regular fiber. The perturbation can be chosen small, so that $[f(N)]$ is unchanged.

Compose $f$ with the projection map $\varphi: M \to \Sigma_g$ of the Seifert fibration. For a regular value $x$ of $\Sigma_g$ whose fiber is transverse to $f$, the preimages of $x$ are the intersection points of $f(N)$ with that fiber, so the mod 2 degree of $\varphi \circ f: N \to \Sigma_g$ is $\psi([f(N)])=\psi(\beta)=1$.

Suppose the singular points of the Seifert fibration on $\Sigma_g$ are $P_i\;(1\leqslant i \leqslant n)$. For each $i$, remove a small disk neighborhood $D(P_i)$ from $\Sigma_g$, remove the vertical solid torus $\varphi^{-1}(D(P_i))$ from the Seifert manifold $M$, and also remove $(\varphi \circ f)^{-1}(D(P_i))$ from the surface $N$. We choose the disks $D(P_i)$ so small that $\varphi \circ f$ is transverse to their boundary circles. This gives a map $\hat g:\hat N\to\hat B$. Here $\hat N=N \setminus \bigcup_{i=1}^n(\varphi \circ f)^{-1}(D(P_i))$ and $\hat B=\Sigma_g \setminus \bigcup_{i=1}^n D(P_i)$. Both surfaces are connected: $\hat B$ is $\Sigma_g$ with the interiors of $n$ disjoint disks removed, and $\hat N$ is a connected closed surface with the interiors of finitely many disjoint disks removed. The map $\hat g$ carries $\partial\hat N$ into $\partial\hat B$. By naturality of the fundamental class, its relative mod 2 degree is the mod 2 degree of $\varphi \circ f$, which is nonzero. Now $\hat B$ has Euler characteristic $2-2g-n$. So Lemma \ref{lem-degree-bound} gives $\chi(\hat N)\leqslant 2-2g-n$.

Next consider the boundary circles of $f(N)$ on the solid torus neighborhoods. A simple closed curve on $\partial N_i$ bounds a disk in $N_i$ precisely when it is a meridian, that is, when its longitude coefficient is $0$. Each boundary circle of a piece removed from $N$ maps onto one of the circles $\partial D(P_i)$ with degree $1$, so its class in $\partial N_i$ is that of the boundary circle of the horizontal surface $\Sigma_g-nD^2$, which runs $-q_i$ times in the longitude direction by the computation in the proof of Proposition \ref{prop-upper-bound-orientable}. Since no $q_i$ vanishes, none of these curves is a meridian, and therefore no piece removed from $N$ is a disk. Hence each of them has Euler characteristic at most $0$. Since the Euler characteristic is additive under gluing along circles, we have $\chi(N)\leqslant\chi(\hat N)\leqslant 2-2g-n$. Therefore, for every element $\beta$ of $H_2(M;\mathbb{Z}_2)\setminus L$, we have $\|\beta\|^{\mathrm{im}}_{\mathbb{Z}_2} \geqslant 2g+n-2$.
\end{proof}

\begin{proof}[Proof of Theorem \ref{thm-main-orientable}]
This follows from Proposition \ref{prop-upper-bound-orientable} and Proposition \ref{prop-lower-bound-orientable}, together with the treatment of the cases $n=0$ and $g=0$ following Proposition \ref{prop-h2-orientable}.
\end{proof}

\section{Seifert Fibered Manifolds with Non-orientable Base}\label{sec-nonorientable}

Let $M$ be a closed Seifert fibered manifold whose base surface is non-orientable and whose total space is orientable. We write $M=(N_g,(p_1,q_1),\cdots,(p_n,q_n))$ with $g\geqslant 1$, where $N_g$ is the connected sum of $g$ projective planes. Since $N_g$ is non-orientable, the product $B_0\times S^1$ of the orientable case has to be replaced by the twisted circle bundle $N_g\overset{\sim}{\times}S^1$ over $N_g$, which is the double of a twisted $I$-bundle over $N_g$. $M$ is obtained from it by removing the $n$ vertical solid tori over small disks $D(P_1),\cdots,D(P_n)$ centred at the singular points $P_1,\cdots,P_n$ of the base orbifold, and gluing them back with gluing maps determined by the pairs $(p_i,q_i)$ as in Section \ref{sec-orientable}. When $n=0$, $M$ is this twisted circle bundle. We denote by $\varphi:M\to N_g$ the projection of the Seifert fibration.

For $M=(N_g,(p_1,q_1),\cdots,(p_n,q_n))$, a straightforward computation gives its fundamental group and its homology groups as follows:
$$
\begin{array}{rcl}
\pi_1(M) = \langle c_1,\cdots,c_g,h_1,\cdots,h_n,v & \mid & vc_iv^{-1}=c_i^{-1}\ (i = 1,\cdots,g),\\
& & h_jv = vh_j,\ h_j^{p_j}v^{q_j} = 1\ (j=1,\cdots,n), \\
& & h_1\cdots h_nc_1^2\cdots c_g^2 = 1 \rangle
\end{array}
$$
$$
\begin{array}{l}
H_1(M;\mathbb{Z}) = \langle [c_1],\cdots,[c_g],[h_1],\cdots,[h_n],[v] \;\; \mid \; 2[c_i]=0\ (i=1,\cdots,g) \\
\quad p_j[h_j] + q_j[v]=0\ (j=1,\cdots,n), \quad [h_1] + \cdots + [h_n]=0 \;\rangle
\end{array}
$$

Below, Proposition \ref{prop-h2-nonorientable} computes $H_2(M;\mathbb{Z}_2)$, Proposition \ref{prop-upper-bound-nonorientable} gives upper bounds for the immersed norm, and Proposition \ref{prop-lower-bound-nonorientable} shows that these bounds are attained.

\begin{Proposition}\label{prop-h2-nonorientable}
Let $M = (N_g,(p_1,q_1),\cdots,(p_n,q_n))$ be a Seifert fibered manifold with $g+n \geqslant 3$ and $n\geqslant 1$. Then:
\begin{enumerate}
\item[(1)] if $p_1,\cdots,p_n$ are all odd, then
\begin{enumerate}
\item[(1.1)] when $q_1+\cdots+q_n$ is odd, we have $H_2(M;\mathbb{Z}_2) \cong \mathbb{Z}_2^{g}$;

\item[(1.2)] when $q_1+\cdots+q_n$ is even, we have $H_2(M;\mathbb{Z}_2) \cong \mathbb{Z}_2^{g+1}$;
\end{enumerate}

\item[(2)] if among $p_1,\cdots,p_n$ there are $k \geqslant 1$ even numbers, then we have $H_2(M;\mathbb{Z}_2) \cong \mathbb{Z}_2^{g+k-1}$.
\end{enumerate}
\end{Proposition}

\begin{proof}
The proof is analogous to that of Proposition \ref{prop-h2-orientable}.
\end{proof}

We use the notation $[v]$, $i_2$ and $\psi$ of Section \ref{sec-orientable}. Recall $\varphi: M \to N_g$ is the projection map. For a simple closed curve $c\subset N_g$ avoiding the singular points, the pre-image $\varphi^{-1}(c)$ is the vertical surface over $c$. It is a Klein bottle of Euler characteristic $0$ when $c$ is one-sided. Let $c_1,\cdots,c_g$ be one-sided simple closed curves, one on each of the projective planes in the connect sum. The mod 2 homology classes of $c_1,\cdots,c_g$ form a basis of $H_1(N_g;\mathbb{Z}_2)$. Define $L\subset H_2(M;\mathbb{Z}_2)$ to be the subgroup generated by the classes of the $g$ vertical Klein bottles $\varphi^{-1}(c_1),\cdots,\varphi^{-1}(c_g)$. Since $i_2(\varphi^{-1}(c_i),[c_j])=\delta_{ij}$, these $g$ classes are independent, so $L\cong\mathbb{Z}_2^{g}$, and every $\alpha\in L$ is represented by a disjoint union of these Klein bottles, so $\|\alpha\|^{\mathrm{im}}_{\mathbb{Z}_2}=0$. Each Klein bottle generating $L$ is disjoint from the fibers over the points of $N_g$ outside the corresponding curve, so $L\subset\ker\psi$.

When $g+n<3$, that is, when $g=1$ and $n=0$, when $g=n=1$, or when $g=2$ and $n=0$, the manifold $M$ admits another Seifert fibration with orientable base (Remark \ref{rem-smallg}) and is covered by Section \ref{sec-orientable}.

When $n=0$ and $g\geqslant 3$, the treatment of the case $n=0$, $g\geqslant 1$ following Proposition \ref{prop-h2-orientable} applies verbatim, with $2g$ replaced by $g$ and $\Sigma_g$ by $N_g$. Here $M$ is the twisted circle bundle, and the presentation above, reduced modulo 2, gives $H_1(M;\mathbb{Z}_2)\cong\mathbb{Z}_2^{g+1}$, so $\dim H_2(M;\mathbb{Z}_2)=g+1$ by Poincar\'e duality. The horizontal surface $N_g$ is the zero section of one of the twisted $I$-bundles whose double is $M$. It has Euler characteristic $2-g$ and meets every regular fiber once, so $i_2([N_g],[v])=1$ and $\psi$ is nonzero, and $L$ has index $2$, so $\ker\psi=L$. Lemma \ref{lem-degree-bound} gives $\chi\leqslant 2-g$ for a surface that maps onto $N_g$ with nonzero mod 2 degree. Thus every $\beta\in H_2(M;\mathbb{Z}_2)\setminus L$ satisfies $\|\beta\|^{\mathrm{im}}_{\mathbb{Z}_2}=g-2=g+n-2$.

We may therefore assume that $M=(N_g,(p_1,q_1),\cdots,(p_n,q_n))$ satisfies $g+n\geqslant 3$ and $n\geqslant 1$.

\begin{Proposition}\label{prop-upper-bound-nonorientable}
Let $M = (N_g,(p_1,q_1),\cdots,(p_n,q_n))$ be a Seifert fibered manifold with $g+n \geqslant 3$ and $n\geqslant 1$. Then:
\begin{enumerate}
\item[(1)] if $p_1,\cdots,p_n$ are all odd, then
\begin{enumerate}
\item[(1.1)] when $q_1+\cdots+q_n$ is odd, every element $\alpha$ of $H_2(M;\mathbb{Z}_2) \cong \mathbb{Z}_2^{g}$ satisfies $\|\alpha\|^{\mathrm{im}}_{\mathbb{Z}_2}=0$;

\item[(1.2)] when $q_1+\cdots+q_n$ is even, every element $\alpha$ of $L$ satisfies $\|\alpha\|^{\mathrm{im}}_{\mathbb{Z}_2}=0$, and every element $\beta$ of $H_2(M;\mathbb{Z}_2)\setminus L$ satisfies $\|\beta\|^{\mathrm{im}}_{\mathbb{Z}_2} \leqslant g+n-2$;
\end{enumerate}

\item[(2)] if among $p_1,\cdots,p_n$ there are $k \geqslant 1$ even numbers, then every element $\alpha\in H_2(M;\mathbb{Z}_2)$ satisfies $\|\alpha\|^{\mathrm{im}}_{\mathbb{Z}_2}=0$.
\end{enumerate}
\end{Proposition}

\begin{proof}
The proof is that of Proposition \ref{prop-upper-bound-orientable} with the following replacements, everything else carrying over verbatim: $\Sigma_g$ becomes $N_g$, the $2g$ curves $a_1,b_1,\cdots,a_g,b_g$ become the $g$ one-sided curves $c_1,\cdots,c_g$, the vertical tori become the vertical Klein bottles over them, and $2g+n$ becomes $g+n$. In case (1.2) the surface $F$ is accordingly built from the horizontal surface $N_g-nD^2$, of Euler characteristic $2-g-n$, and the $n$ immersed M\"obius bands, hence has genus $g+n$. Case (2) needs no change, since its construction takes place inside the solid tori.
\end{proof}

Thus it remains only to show that the inequality in Proposition \ref{prop-upper-bound-nonorientable} (1.2) is in fact an equality. Lemma \ref{lem-degree-bound} applies without change, since it is stated for arbitrary compact surfaces and its mod 2 formulation needs no orientation.

\begin{Proposition}\label{prop-lower-bound-nonorientable}
Let $M = (N_g,(p_1,q_1),\cdots,(p_n,q_n))$ be a Seifert fibered manifold with $g+n \geqslant 3$ and $n\geqslant 1$. Suppose that all $p_i$ are odd, that no $q_i$ vanishes, and that $q_1+\cdots+q_n$ is even, so that $H_2(M;\mathbb{Z}_2) \cong \mathbb{Z}_2^{g+1} = L \oplus \mathbb{Z}_2$ with $L\cong \mathbb{Z}_2^{g}$. Then every element $\beta$ of $H_2(M;\mathbb{Z}_2)\setminus L$ satisfies $\|\beta\|^{\mathrm{im}}_{\mathbb{Z}_2} \geqslant g+n-2$.
\end{Proposition}

\begin{proof}
The proof is that of Proposition \ref{prop-lower-bound-orientable}, with $N_g$ in place of $\Sigma_g$ and $2-g-n$ in place of $2-2g-n$, everything else carrying over verbatim. The one point to note is that $N_g$ is non-orientable, so the map $\varphi \circ f$ has no integral degree and only its mod 2 degree is available. This is enough: for a regular value $x$ of $N_g$ whose fiber is transverse to $f$, the preimages of $x$ are exactly the intersection points of $f(N)$ with that fiber, so that mod 2 degree is $\psi(\beta)=1$, and the estimate of $\chi(N)$ comes from Lemma \ref{lem-degree-bound}, which is stated in mod 2 homology and needs no orientation.
\end{proof}

\begin{proof}[Proof of Theorem \ref{thm-main-nonorientable}]
This follows from Proposition \ref{prop-upper-bound-nonorientable} and Proposition \ref{prop-lower-bound-nonorientable}, together with the treatment of the case $n=0$ following Proposition \ref{prop-h2-nonorientable}.
\end{proof}

\section{The Case of Torus Bundles}\label{sec-torusbundle}

Suppose that $M = T^2 \times [0,1] / ((x,0) \sim (\varphi(x),1))$ is a torus bundle, and that the matrix of the induced homomorphism $\varphi_*: H_1(T^2) \to H_1(T^2)$ is $\left(\begin{smallmatrix}a & b \\ c & d\end{smallmatrix}\right)$. As explained in the introduction, we only need to consider the case where $\varphi$ is an Anosov map. The mod 2 residue class $\left[\left(\begin{smallmatrix}a & b \\ c & d\end{smallmatrix}\right)\right]_2$ is then one of the following six matrices, which together exhaust $SL(2;\mathbb{Z}_2)$: $\left(\begin{smallmatrix}1 & 1 \\ 1 & 0\end{smallmatrix}\right)$, $\left(\begin{smallmatrix}0 & 1 \\ 1 & 1\end{smallmatrix}\right)$, $\left(\begin{smallmatrix}1 & 1 \\ 0 & 1\end{smallmatrix}\right)$, $\left(\begin{smallmatrix}1 & 0 \\ 1 & 1\end{smallmatrix}\right)$, $\left(\begin{smallmatrix}0 & 1 \\ 1 & 0\end{smallmatrix}\right)$, $\left(\begin{smallmatrix}1 & 0 \\ 0 & 1\end{smallmatrix}\right)$.

A straightforward computation gives the following result about the homology groups:

\begin{Proposition}\label{prop-h2-torusbundle}
Let $M = T^2 \times [0,1] / ((x,0) \sim (\varphi(x),1))$ be the torus bundle given by an Anosov map $\varphi$. For such a bundle we have $b_1(M)=1$, and Poincar\'e duality gives $H_2(M;\mathbb{Z}) \cong H^1(M;\mathbb{Z}) \cong \mathbb{Z}$. For the mod 2 homology groups we have:
\begin{enumerate}
\item[(1)] when $\left[\left(\begin{smallmatrix}a & b \\ c & d\end{smallmatrix}\right)\right]_2=\left(\begin{smallmatrix}1 & 1 \\ 1 & 0\end{smallmatrix}\right)$ or $\left(\begin{smallmatrix}0 & 1 \\ 1 & 1\end{smallmatrix}\right)$, we have $H_2(M;\mathbb{Z}_2)\cong\mathbb{Z}_2$;

\item[(2)] when $\left[\left(\begin{smallmatrix}a & b \\ c & d\end{smallmatrix}\right)\right]_2=\left(\begin{smallmatrix}1 & 1 \\ 0 & 1\end{smallmatrix}\right)$, $\left(\begin{smallmatrix}1 & 0 \\ 1 & 1\end{smallmatrix}\right)$ or $\left(\begin{smallmatrix}0 & 1 \\ 1 & 0\end{smallmatrix}\right)$, we have $H_2(M;\mathbb{Z}_2)\cong\mathbb{Z}_2^2$;

\item[(3)] when $\left[\left(\begin{smallmatrix}a & b \\ c & d\end{smallmatrix}\right)\right]_2=\left(\begin{smallmatrix}1 & 0 \\ 0 & 1\end{smallmatrix}\right)$, we have $H_2(M;\mathbb{Z}_2)\cong\mathbb{Z}_2^3$.
\end{enumerate}
\end{Proposition}

We split the proof of Theorem \ref{thm-torusbundle} into the following propositions.

\begin{Proposition}\label{prop-torusbundle-upper}
Let $M$ be the torus bundle over the circle arising from an Anosov map. The unique nonzero element $\alpha$ of $H_2(M;\mathbb{Z}) \otimes \mathbb{Z}_2 \cong \mathbb{Z}_2$ satisfies $\|\alpha\|^{\mathrm{im}}_{\mathbb{Z}_2}=0$, and all the remaining nonzero elements $\beta$ of $H_2(M;\mathbb{Z}_2)$ (if there are any) satisfy $\|\beta\|^{\mathrm{im}}_{\mathbb{Z}_2} \leqslant 1$.
\end{Proposition}

\begin{proof}
We only need to construct surfaces representing the nonzero homology classes in $H_2(M;\mathbb{Z}_2)$ and compute their Euler characteristics.

First, the fiber torus represents the homology class $\alpha\in H_2(M;\mathbb{Z}) \otimes \mathbb{Z}_2$, so that $\|\alpha\|^{\mathrm{im}}_{\mathbb{Z}_2}=0$.

For a primitive class $\gamma\in H_1(T^2;\mathbb{Z})$ the mod 2 intersection number of $\gamma$ with $\varphi(\gamma)$ equals $\det(\gamma,\varphi(\gamma))$. In the paragraphs below we choose $\gamma$ so that $\varphi(\gamma)\equiv \gamma$ in $H_1(T^2;\mathbb{Z}_2)$, so that this number is even, while $\varphi(\gamma)\neq \pm\gamma$ because $\varphi$ is an Anosov map.

When $\left[\left(\begin{smallmatrix}a & b \\ c & d\end{smallmatrix}\right)\right]_2=\left(\begin{smallmatrix}1 & 1 \\ 0 & 1\end{smallmatrix}\right)$, the class $\gamma=\left(\begin{smallmatrix}1\\0\end{smallmatrix}\right)$ satisfies $\varphi(\gamma)\equiv \gamma$, so the mod 2 intersection number of $\gamma$ with its image $\varphi(\gamma)$ is even, and $\varphi(\gamma)\neq \pm\gamma$. By Proposition \ref{prop-punctured-band} we can construct in $T^2 \times [0,1]$ an immersed punctured M\"obius band whose two boundary curves are, respectively, the curve $(\gamma,0)$ on $T^2 \times 0$ and the curve $(\varphi(\gamma),1)$ on $T^2 \times 1$. This gives an immersed non-orientable surface $i: \mathbb{R}P^2 \# \mathbb{R}P^2 \# \mathbb{R}P^2 \looparrowright M$ of Euler characteristic $-1$. It intersects the fiber torus in a non-trivial curve, hence it represents a mod 2 homology class $\beta$ lying outside $H_2(M;\mathbb{Z})\otimes \mathbb{Z}_2$, and we have $\|\beta\|^{\mathrm{im}}_{\mathbb{Z}_2} \leqslant 1$. Taking the union of the above surface with the fiber torus gives a representing surface of $\beta + \alpha$, so that $\|\beta + \alpha\|^{\mathrm{im}}_{\mathbb{Z}_2} \leqslant 1$. By the same intersection-number argument, $\alpha,\beta,\alpha+\beta$ represent the three nonzero homology classes in $H_2(M;\mathbb{Z}_2)\cong \mathbb{Z}_2^2$ in this case.

When $\left[\left(\begin{smallmatrix}a & b \\ c & d\end{smallmatrix}\right)\right]_2=\left(\begin{smallmatrix}1 & 0 \\ 1 & 1\end{smallmatrix}\right)$, it is the class $\gamma=\left(\begin{smallmatrix}0\\1\end{smallmatrix}\right)$ that satisfies $\varphi(\gamma)\equiv \gamma$. When $\left[\left(\begin{smallmatrix}a & b \\ c & d\end{smallmatrix}\right)\right]_2=\left(\begin{smallmatrix}0 & 1 \\ 1 & 0\end{smallmatrix}\right)$, it is the class $\gamma=\left(\begin{smallmatrix}1\\1\end{smallmatrix}\right)$ that satisfies $\varphi(\gamma)\equiv \gamma$. These cases are analogous to the previous one.

When $\left[\left(\begin{smallmatrix}a & b \\ c & d\end{smallmatrix}\right)\right]_2=\left(\begin{smallmatrix}1 & 0 \\ 0 & 1\end{smallmatrix}\right)$, the classes $\left(\begin{smallmatrix}1\\0\end{smallmatrix}\right)$, $\left(\begin{smallmatrix}0\\1\end{smallmatrix}\right)$ and $\left(\begin{smallmatrix}1\\1\end{smallmatrix}\right)$ all satisfy $\varphi(\gamma)\equiv \gamma$. Hence we can construct 3 immersed surfaces in the same way. These surfaces, together with their unions with the fiber torus, give 6 immersed surfaces, each of Euler characteristic $-1$. By computing intersection numbers with curves on the fiber torus, one sees that these surfaces represent distinct mod 2 homology classes, and that they are different from the mod 2 homology class represented by the fiber torus. Hence every mod 2 homology class $\beta$ represented by them satisfies $\|\beta\|^{\mathrm{im}}_{\mathbb{Z}_2} \leqslant 1$.
\end{proof}

\begin{Proposition}\label{prop-torusbundle-lower}
Let $M$ be the torus bundle over the circle arising from an Anosov map, and let $L = H_2(M;\mathbb{Z}) \otimes \mathbb{Z}_2$. Then every element $\beta$ of $H_2(M;\mathbb{Z}_2)\setminus L$ satisfies $\|\beta\|^{\mathrm{im}}_{\mathbb{Z}_2} \geqslant 1$.
\end{Proposition}

\begin{proof}
Since $\beta\notin L=H_2(M;\mathbb{Z})\otimes\mathbb{Z}_2$, the vertical component of $\beta$ is nonzero, that is, every closed surface representing $\beta$ meets the fiber torus in an odd number of curves.

For such a surface $F$, cutting $M$ along the fiber torus yields an immersed surface $i : F_1 \looparrowright T^2 \times [0,1]$, and $F$ is recovered from $F_1$ by gluing its boundary circles in pairs, so that $\chi(F)=\chi(F_1)$. We first discard the components of $F_1$ that contribute nothing. A component without boundary circle is a closed surface immersed in $T^2 \times [0,1]$: if it is a sphere or a projective plane it can be eliminated by a homotopy, since $\pi_2(T^2 \times [0,1])=0$ and $\pi_1(T^2 \times [0,1])$ is torsion-free, and if it is a torus or a Klein bottle its Euler characteristic is $0$, so it may be ignored. If a component of $F_1$ has only one boundary circle, then that circle is null-homologous in the boundary torus and hence null-homotopic there, because $\pi_1$ of the boundary torus injects into $\pi_1(T^2 \times [0,1])$, so the immersion of that component can again be eliminated by a homotopy.

Hence we may assume that every connected component of $F_1$ has at least 2 boundary circles. Then the Euler characteristic of each component of $F_1$ is at most $0$, and if it equals $0$ then that component is an annulus. Moreover, an annulus component has two boundary circles that are homologous in $T^2$, and these two circles are identified with each other by the gluing map $\varphi$. Such a component would therefore give a nonzero class $c\in H_1(T^2;\mathbb{Z})$ with $\varphi(c)=\pm c$, which is impossible since $\varphi$ is an Anosov map. It follows that at least one component of $F_1$ is not an annulus, so such a component has $\chi\leqslant -1$. Hence $\|\beta\|^{\mathrm{im}}_{\mathbb{Z}_2} \geqslant 1$.
\end{proof}

\section{The Case of Torus Semi-bundles}\label{sec-semibundle}

We regard a torus semi-bundle as the space obtained by the following gluing. Take two copies of the Seifert fibered manifold $(D^2,(2,1),(2,-1))$ (homeomorphic to $K \overset{\sim}{\times} I$, where $K$ is the Klein bottle) and one copy of the thickened torus $T^2 \times [0,1]$. On the boundary torus of the $i$-th copy of $(D^2,(2,1),(2,-1))$ ($i = 1,2$), take a homology basis $\left\{[l_i]=\left(\begin{smallmatrix}1\\0\end{smallmatrix}\right),[m_i]=\left(\begin{smallmatrix}0\\1\end{smallmatrix}\right)\right\}$, where $[l_i]$ is represented by a regular fiber on the boundary torus and $[m_i]$ is represented by a curve $\partial D^2\times \{\mathrm{pt}\}$. For the homology group of the thickened torus, take a basis $(l_3,m_3)$.

The gluing of the boundary of the first copy of $(D^2,(2,1),(2,-1))$ to $T^2 \times 0$ satisfies $[l_1] \mapsto [l_3]$, $[m_1] \mapsto [m_3]$. The gluing $\varphi$ from $T^2 \times 1$ to the boundary of the second copy of $(D^2,(2,1),(2,-1))$ satisfies $[l_3] \mapsto p [l_2] + q [m_2]$, $[m_3] \mapsto r[l_2] + s [m_2]$. Denote by $[\varphi]_2=\left[\left(\begin{smallmatrix} p & r \\ q & s\end{smallmatrix}\right)\right]_2$ the mod 2 residue class of the matrix $\left(\begin{smallmatrix} p & r \\ q & s\end{smallmatrix}\right)$ of $\varphi$.

We may also regard the thickened torus $T^2 \times [0,1]$ as a regular neighborhood of the boundary of the first copy, and hence regard the torus semi-bundle as two copies of $(D^2,(2,1),(2,-1))$ glued along their boundaries by $\varphi$.

\begin{Proposition}\label{prop-h2-semibundle}
Let $M$ be the torus semi-bundle obtained as above, where $\varphi$ is an Anosov map. For integral homology, $H_2(M;\mathbb{Z})=0$. For mod 2 homology we have:
\begin{enumerate}
\item[(1)] when $\left[\left(\begin{smallmatrix} p & r \\ q & s\end{smallmatrix}\right)\right]_2 = \left(\begin{smallmatrix}1 & 0 \\ 0 & 1\end{smallmatrix}\right)$ or $\left(\begin{smallmatrix}1 & 1 \\ 0 & 1\end{smallmatrix}\right)$, $H_2(M;\mathbb{Z}_2) \cong \mathbb{Z}_2^3$;

\item[(2)] when $\left[\left(\begin{smallmatrix} p & r \\ q & s\end{smallmatrix}\right)\right]_2 =\left(\begin{smallmatrix}1 & 1 \\ 1 & 0\end{smallmatrix}\right)$, $\left(\begin{smallmatrix}1 & 0 \\ 1 & 1\end{smallmatrix}\right)$, $\left(\begin{smallmatrix}0 & 1 \\ 1 & 1\end{smallmatrix}\right)$ or $\left(\begin{smallmatrix}0 & 1 \\ 1 & 0\end{smallmatrix}\right)$, $H_2(M;\mathbb{Z}_2) \cong \mathbb{Z}_2^2$.
\end{enumerate}
\end{Proposition}

\begin{proof}
By Poincar\'e duality and the universal coefficient theorem, we only need to compute $H_1(M;\mathbb{Z})$. Since $\varphi$ is an Anosov map, the entry $r$ of the gluing matrix $\left(\begin{smallmatrix} p & r \\ q & s\end{smallmatrix}\right)$ satisfies $r \neq 0$: indeed, $r=0$ together with $\det=ps-qr=\pm 1$ would give $p, s = \pm 1$, so the matrix is triangular with eigenvalues $p$ and $s$ both of modulus $1$, and $\varphi$ would not be Anosov.

Regard the torus semi-bundle as two copies of $(D^2,(2,1),(2,-1))\cong K \overset{\sim}{\times} I$ glued along their boundaries by $\varphi$. The fundamental group of both the Klein bottle $K$ and its twisted $I$-bundle $K \overset{\sim}{\times} I$ is $\langle a,b \mid a^2=b^2 \rangle$. The boundary of $K \overset{\sim}{\times} I$ is a torus, whose fundamental group can be generated by $a^2$ and $b^{-1}a$. The regular fiber of the Seifert fibration $(D^2,(2,1),(2,-1))$ represents the element $a^2$ (which equals $b^2$). $\partial D^2$ represents the element $b^{-1}a$. To verify $a^2$ and $b^{-1}a$ commute, we only need to note that
$$
\begin{array}{rl}
a^2(b^{-1}a)(a^2)^{-1}(b^{-1}a)^{-1} & = aa(b^{-1}a)a^{-1}a^{-1}(a^{-1}b) \\
& = aab^{-1}a^{-1}a^{-1}b=(bb)b^{-1}(b^{-1}b^{-1})b=1.
\end{array}
$$
When $M$ is the result of two copies of $(D^2,(2,1),(2,-1))$ glued by $\varphi$, a presentation of $H_1(M;\mathbb{Z})$ is
$$
\begin{array}{rl}
\langle [a_1],[a_2],[b_1],[b_2] \mid & 2[a_1]=2[b_1], \;2[a_2]=2[b_2], \\
& 2[a_1] = p (2[a_2]) + q (-[b_2]+[a_2]), \\
& -[b_1]+[a_1] = r(2[a_2]) + s (-[b_2]+[a_2])\rangle
\end{array}
$$

A straightforward computation gives $H_1(M;\mathbb{Z}) \cong \mathbb{Z}_{4|r|}\oplus \mathbb{Z}_2 \oplus \mathbb{Z}_2$ when $q$ is even (hence $H_2(M;\mathbb{Z})=0$ and $\operatorname{Tor}(H_1(M;\mathbb{Z}),\mathbb{Z}_2) \cong \mathbb{Z}_2^3$), and $H_1(M;\mathbb{Z})\cong \mathbb{Z}_{4|r|}\oplus\mathbb{Z}_4$ when $q$ is odd (hence $H_2(M;\mathbb{Z})=0$ and $\operatorname{Tor}(H_1(M;\mathbb{Z}),\mathbb{Z}_2) \cong \mathbb{Z}_2^2$).
\end{proof}

We use the notation $i_2$ of Section \ref{sec-orientable}. In the two copies of $(D^2,(2,1),(2,-1)) \cong K \overset{\sim}{\times} I$, the core Klein bottles $K_1,K_2$ are embedded and disjoint in $M$. We define $L\subset H_2(M;\mathbb{Z}_2)$ to be the subgroup generated by the classes $[K_1],[K_2]$. A closed curve $\gamma_i$ in the $i$-th copy that crosses $K_i$ once, for example the union of an $I$-fiber of the twisted bundle with a path on the boundary torus, satisfies $i_2([K_j],[\gamma_i])=\delta_{ij}$, so $[K_1]$ and $[K_2]$ are independent and $L\cong\mathbb{Z}_2^2$. Obviously every $\alpha\in L$ satisfies $\|\alpha\|^{\mathrm{im}}_{\mathbb{Z}_2}=0$.

We split the proof of Theorem \ref{thm-semibundle} into the following propositions.

\begin{Proposition}\label{prop-semibundle-upper}
Let $M$ be the torus semi-bundle obtained as above, where $\varphi$ is an Anosov map. Then every element $\beta$ of $H_2(M;\mathbb{Z}_2)\setminus L$ (if there is any) satisfies $\|\beta\|^{\mathrm{im}}_{\mathbb{Z}_2} \leqslant 1$.
\end{Proposition}

\begin{proof}
It suffices to find an immersed representing surface for a given element $\beta$ of $H_2(M;\mathbb{Z}_2)\setminus L$.

When $q$ in the matrix $\left(\begin{smallmatrix} p & r \\ q & s\end{smallmatrix}\right)$ describing $\varphi$ is odd, Proposition \ref{prop-h2-semibundle} gives $H_2(M;\mathbb{Z}_2)\cong \mathbb{Z}_2^2=L$, and there is nothing left to prove.

When $q$ in the matrix $\left(\begin{smallmatrix} p & r \\ q & s\end{smallmatrix}\right)$ describing $\varphi$ is even, Proposition \ref{prop-h2-semibundle} gives $H_2(M;\mathbb{Z}_2)\cong \mathbb{Z}_2^3$. In each of the two copies of $(D^2,(2,1),(2,-1))$ take a vertical regular fiber lying on the boundary. These two fibers bound embedded M\"obius bands in the two copies of $(D^2,(2,1),(2,-1))$ respectively. In the basis $(l_3,m_3)$ of $H_1(T^2)$, the first of these two curves is $[l_3]$, and the second is $\pm (s[l_3]-q[m_3])$, obtained by inverting the gluing matrix. Hence the two curves are not $\pm$ each other, and their mod 2 intersection number is an even number $q$. By Proposition \ref{prop-punctured-band}, there is an immersed punctured M\"obius band in $T^2 \times [0,1]$ bounded by these two curves. This allows us to assemble a surface $F \cong \mathbb{R}P^2\#\mathbb{R}P^2\#\mathbb{R}P^2 \looparrowright M$. By computing intersection numbers one sees that $[F] \not\in L$. Taking unions of $F$ with representing surfaces of the various elements of $L$ then gives immersed representing surfaces of all the elements of $H_2(M;\mathbb{Z}_2)$, and every element $\beta$ of $H_2(M;\mathbb{Z}_2) \setminus L$ satisfies $\|\beta\|^{\mathrm{im}}_{\mathbb{Z}_2}\leqslant 1$.
\end{proof}

\begin{Proposition}\label{prop-semibundle-lower}
Let $M$ be the torus semi-bundle obtained as above, where $\varphi$ is an Anosov map. When $q$ in the matrix $\left(\begin{smallmatrix} p & r \\ q & s\end{smallmatrix}\right)$ describing $\varphi$ is even, every element $\beta$ of $H_2(M;\mathbb{Z}_2) \setminus L$ satisfies $\|\beta\|^{\mathrm{im}}_{\mathbb{Z}_2}\geqslant 1$.
\end{Proposition}

\begin{proof}
Let an immersed surface $i:F \looparrowright M$ satisfy $[i(F)]=\beta$. Then $i(F)$ necessarily intersects the fiber torus. In fact, if $i(F)$ were disjoint from the fiber torus, then it would lie in one copy of $(D^2,(2,1),(2,-1))\cong K \overset{\sim}{\times}I$. This $I$-bundle deformation retracts onto its core Klein bottle representing an element in $L$. So $\beta=[i(F)]$ belongs to $L$, contradicting $\beta\in H_2(M;\mathbb{Z}_2)\setminus L$.

Cutting $M$ open into two copies of $(D^2,(2,1),(2,-1)) \cong K \overset{\sim}{\times}I$, we obtain two immersed surfaces $i_1:F_1 \looparrowright K \overset{\sim}{\times}I$ and $i_2:F_2 \looparrowright K \overset{\sim}{\times}I$ (each $F_i$ can be disconnected), and $F_1,F_2$ can be glued together to recover $F$, so $\chi(F)=\chi(F_1) + \chi(F_2)$.

We now estimate $\chi(F_1)$ and $\chi(F_2)$. First we exclude the possibility that $F_1$ or $F_2$ has a connected component of positive Euler characteristic. In fact, such a component can only be a sphere, a projective plane or a disk. If the component is a projective plane, then we claim it represents a zero class and hence can be deleted. Indeed, if it represented a nonzero class in $H_2(M;\mathbb{Z}_2)$, then, since $H_2(M;\mathbb{Z})=0$ by Proposition \ref{prop-h2-semibundle}, that class would lie outside $H_2(M;\mathbb{Z})\otimes\mathbb{Z}_2$, and Proposition \ref{prop-pi2-detection} would give an element of order $2$ in $\pi_1(M)$. However $\pi_1(M)$ is torsion-free, a contradiction. If the component is a sphere, then $\pi_2(K \overset{\sim}{\times}I)=0$ implies that its image is contractible in $K \overset{\sim}{\times}I$, so it contributes nothing to the mod 2 homology class. If the component is a disk, then we claim that it can be homotoped off $K \overset{\sim}{\times}I$. Indeed, the boundary curve of this disk is null-homotopic in $K \overset{\sim}{\times}I$, so its image in $\pi_1(K \overset{\sim}{\times}I)\cong \pi_1(K)$ is trivial. Since $\partial(K \overset{\sim}{\times}I)$ is a two-fold cover of $K$, the map $\pi_1(\partial(K \overset{\sim}{\times}I)) \to \pi_1(K \overset{\sim}{\times}I)$ is injective. Hence the curve is null-homotopic on the torus $\partial(K \overset{\sim}{\times}I)$ itself, and there is an immersed disk lying entirely in $\partial(K \overset{\sim}{\times}I)$ with this immersed curve as boundary. The two immersed disks can be glued into an immersed sphere, which is contractible since $\pi_2(K \overset{\sim}{\times}I)$ is trivial. After contracting this immersed sphere to a point, we can push the disk component off $K \overset{\sim}{\times}I$. Thus $\chi(F_1)\leqslant 0$ and $\chi(F_2)\leqslant 0$.

In order to prove $\chi(F_1)+\chi(F_2)\leqslant -1$, it remains to show that $\chi(F_1)=0$ and $\chi(F_2)=0$ cannot hold simultaneously. When this happens, every component of $F_1$ and of $F_2$ has Euler characteristic $0$, so each of them is a torus, a Klein bottle, an annulus or a M\"obius band. The closed ones may be ignored below: having no boundary circle, they play no part in the gluing that closes up $F$. It therefore suffices to discuss the components with boundary circles. For an immersed annulus component, its two boundary curves are homotopic on the torus $\partial(K \overset{\sim}{\times}I)$. Hence we can modify the immersed closed surface $F$ as follows: replace this immersed annulus by an immersed annulus lying on the boundary torus $\partial(K \overset{\sim}{\times}I)$, and then push it slightly outwards. The modified immersed closed surface $F'$ still represents the same homology class as $F$, its Euler characteristic is unchanged, and after cutting it open there is no longer such an immersed annulus component. As for the immersed M\"obius band components, the boundary curve of an immersed M\"obius band in $K \overset{\sim}{\times}I \cong (D^2,(2,1),(2,-1))$ is, in $\pi_1$, the square of its core curve, and since $a^2$ is central in $\pi_1(K)=\langle a,b \mid a^2=b^2\rangle$, the square of any element is homotopic to a power of the regular fiber. If two immersed M\"obius bands in the two copies of $K \overset{\sim}{\times}I$ could be glued along their boundaries, then $\varphi$ would carry one power of the regular fiber to another, hence would preserve the direction of the fiber on $T^2$, so that $q=0$. This contradicts the fact that $\varphi$ is an Anosov map, since $q=0$ together with $\det=\pm1$ gives eigenvalues $\pm1$.
\end{proof}

\section{Immersed Norm and Minimal Genus}\label{sec-imnorm}

For $\alpha \in H_2(M;\mathbb{Z}_2)\setminus (H_2(M;\mathbb{Z})\otimes \mathbb{Z}_2)$, if we restrict the representing surfaces to non-orientable surfaces $N_k$, we can also define the minimal immersed genus
$$
g^{\mathrm{im}}_N(\alpha) = \min\{ k: N_k \looparrowright M \text{ and } [N_k]=\alpha\}.
$$

\begin{Remark}\label{rem-crosscap}
When $M$ is the compact manifold obtained from $S^3$ by removing an open regular neighborhood of a knot and $\alpha$ is a mod 2 homology class relative to the boundary, the quantity $\|\alpha\|^{\mathrm{em}}_{\mathbb{Z}_2}$ is, in most cases, a variant of what is called the ``crosscap number'' or the ``non-orientable genus'' of that knot in the literature. The immersed analogue is defined similarly. \cite{HughesKim} considered the case in which such relative mod 2 homology classes satisfy $\|\alpha\|^{\mathrm{im}}_{\mathbb{Z}_2}=0$ while the immersed crosscap number equals 1.
\end{Remark}

By the definition of the immersed norm, when $\|\alpha\|^{\mathrm{im}}_{\mathbb{Z}_2}=0$ and the simplest representing surface $S$ is connected, there are still four possibilities for $S$: sphere $S^2$, projective plane $\mathbb{R}P^2$, torus $T^2$, and Klein bottle $K$. If $S$ is $S^2$ or $T^2$, then $\alpha \in H_2(M;\mathbb{Z}) \otimes \mathbb{Z}_2$. If $\alpha \in H_2(M;\mathbb{Z}_2) \setminus H_2(M;\mathbb{Z}) \otimes \mathbb{Z}_2$ and $\|\alpha\|^{\mathrm{im}}_{\mathbb{Z}_2}=0$, we need to determine whether $g_N^{\mathrm{im}}(\alpha)=1$ or $g_N^{\mathrm{im}}(\alpha)=2$.

\begin{Proposition}\label{prop-projective-plane-cases}
Let $M$ be a closed, orientable, non-hyperbolic geometric 3-manifold. Then a class $\alpha \in H_2(M;\mathbb{Z}_2) \setminus (H_2(M;\mathbb{Z}) \otimes \mathbb{Z}_2)$ satisfies $g^{\mathrm{im}}_N(\alpha)=1$ if and only if $M$ is either a lens space $L(2k,q)$ or the Seifert manifold $(S^2,(2,1),(2,1),(n,q))$ with $n$ odd.
\end{Proposition}

\begin{proof}
By Proposition \ref{prop-pi2-detection}, there is a class $\alpha \in H_2(M;\mathbb{Z}_2) \setminus H_2(M;\mathbb{Z})\otimes\mathbb{Z}_2$ satisfying $g^{\mathrm{im}}_N(\alpha)=1$ if and only if $\pi_1(M)$ contains an element of order $2$ that is mapped to an element of order $2$ in $H_1(M;\mathbb{Z})$. If $M$ were not elliptic, then $\pi_1(M)$ would be torsion-free by the geometrization theorem, so $M$ must be elliptic.

Next, $M = S^3/G$ with $G = \pi_1(M)$ finite acting freely on $S^3$. A free involution of $S^3$ is the antipodal map, so $G$ has a unique element of order $2$, namely $-1$, and the question is whether $-1$ has trivial image in $H_1(M;\mathbb{Z}) = G^{\mathrm{ab}}$. For the lens space $L(2k,q)$ the group $\pi_1(M)$ is cyclic of order $2k$ and $-1$ has image of order $2$, while for $L(p,q)$ with $p$ odd there is no element of order $2$ at all.

It remains to consider the elliptic manifolds listed in \cite[Table 12.5]{Martelli}, whose invariants are $(S^2,(p_1,q_1),(p_2,q_2),(p_3,q_3))$ with $(p_1,p_2,p_3)$ equal to $(2,2,n)$, $(2,3,3)$, $(2,3,4)$ or $(2,3,5)$. We use the notation of that table, where $C_m$ is the cyclic group of order $m$ and $T^*_{24}$, $O^*_{48}$, $I^*_{120}$, $D^*_{4n}$ are the binary tetrahedral, octahedral, icosahedral and dihedral groups. For $(2,3,3)$ the group is $T^*_{24}\times C_{|5+2q|}$, or a subgroup of index $3$ of $T^*_{24}\times C_{6(5+2q)}$ when $3\mid 5+2q>0$, and in both cases $|H_1| = 3|5+2q|$ is odd, so $-1$ has trivial image. For $(2,3,5)$ the group is $I^*_{120}\times C_{|25+6q|}$ and $|H_1| = |25+6q|$ is odd, so $-1$ again has trivial image. For $(2,3,4)$ the group is $O^*_{48}\times C_{|10+3q|}$, which contains the quaternion subgroup $\{\pm1,\pm i,\pm j,\pm k\}$, in which $-1=[i,j]$, so $-1$ is a commutator and has trivial image.

For $(2,2,n)$ the group is $D^*_{4n}\times C_{n+q}$ when $n+q$ is odd, and a subgroup of index $2$ of $D^*_{4n}\times C_{4(n+q)}$ when $n+q$ is even, which forces $n$ odd. In the former case, we claim $-1$ has trivial image exactly when $n$ is even. In fact, for even $n$, $D^*_{4n}$ contains the quaternion subgroup $\{\pm1,\pm i,\pm j,\pm k\}$ and $-1=[i,j]$. For odd $n$, the abelianization of $D^*_{4n}$ is $\mathbb{Z}_4$ and $-1$, being the square of a generator, has image of order $2$. Multiplying with a cyclic group does not change this, since the commutator subgroup of a product is the product of the commutator subgroups. Neither does passing to a subgroup of finite index, since a commutator of a subgroup is a commutator of the ambient group. Hence, by Proposition \ref{prop-pi2-detection}, such a class exists exactly in the two cases listed in Proposition \ref{prop-projective-plane-cases}.
\end{proof}

Combining the results of Sections \ref{sec-orientable}, \ref{sec-nonorientable}, \ref{sec-torusbundle} and \ref{sec-semibundle}, we in fact also obtain:

\begin{Proposition}\label{prop-genus-two}
Let $M$ be a closed, orientable, non-hyperbolic geometric 3-manifold. Then there exists $\alpha \in H_2(M;\mathbb{Z}_2) \setminus H_2(M;\mathbb{Z}) \otimes \mathbb{Z}_2$ satisfying $g^{\mathrm{im}}_N(\alpha)=2$ if and only if $M$ is of one of the following two types: (1) a Seifert manifold whose base orbifold has at least two singular points of even index, or (2) a torus semi-bundle.
\end{Proposition}

For a class in $H_2(M;\mathbb{Z})$ that can be represented by a connected surface, the simplest such representing surface is necessarily $\pi_1$-injective. By contrast, Propositions \ref{prop-projective-plane-cases} and \ref{prop-genus-two} show that there are manifolds $M$ for which the simplest immersed connected representing surface of an element of $H_2(M;\mathbb{Z}_2)$ is not $\pi_1$-injective.

\begin{Example}\label{eg-noninjective}
For the Seifert manifold $M=(S^2,(2,1),(2,1),(n,q))$ with $n$ even, or $M=(S^2,(2,1),(3,1),(4,-3))$, every nonzero class $\beta \in H_2(M;\mathbb{Z}_2)$ satisfies $\|\beta\|^{\mathrm{im}}_{\mathbb{Z}_2}=0$ and $g_N^{\mathrm{im}}(\beta)=2$, the value $1$ being excluded by Proposition \ref{prop-projective-plane-cases}. Hence any simplest non-orientable immersed surface $S$ representing $\beta$ is homeomorphic to $\mathbb{R}P^2 \# \mathbb{R}P^2$, and $\pi_1(S)$ is infinite. But $\pi_1(M)$ is finite, so $S$ is not $\pi_1$-injective.
\end{Example}

\begin{Remark}\label{rem-embedded-injective}
For a primitive integral homology class, if its simplest embedded representing surface is connected, then that surface is $\pi_1$-injective. For a 2-dimensional mod 2 homology class coming from $\operatorname{Tor}(H_1(M;\mathbb{Z}),\mathbb{Z}_2)$, if its simplest embedded representing surface is connected, it is shown in \cite{BW} that such a surface need not be $\pi_1$-injective. The above example shows that even if we replace embedded representing surfaces by immersed ones, these simplest immersed surfaces may still fail to be $\pi_1$-injective.
\end{Remark}

\textbf{AI declaration.} The authors used Doubao to revise the English for clarity, improve formatting of the manuscript, and find typos. The authors verified the final proofs and takes full responsibility for the content of the paper.

\end{document}